\documentclass[draft]{agujournal2019}
\usepackage{url} 
\usepackage{soul}
\usepackage{graphicx}
\usepackage{amsmath}
\usepackage{amsthm}
\usepackage{algorithm}
\usepackage{algpseudocode} 

\usepackage{cancel}
\usepackage{array}
\usepackage{float}
\usepackage{soul}
\usepackage{breqn}
\usepackage{multirow}

\draftfalse

\journalname{Water Resources Research}

\begin{document}

%
%


\title{Bound-Preserving Flux-Corrected Transport Methods for Solving Richards' Equation}
%
%




\authors{Arnob Barua\affil{1}, Christopher E. Kees\affil{1,2}, Dmitri Kuzmin\affil{3}}


\affiliation{1}{Department of Civil and Environmental Engineering, Louisiana State University}
\affiliation{2}{Center for Computation and Technology, Louisiana State University}
\affiliation{3}{Department of Mathematics, TU Dortmund University}




\correspondingauthor{Christopher E. Kees}{cekees@lsu.edu}



\begin{keypoints}
\item A first-order method for Richards' equation preserving physical bounds for general unstructured meshes and anisotropy is presented.
\item A novel second-order flux-corrected transport method is devised for Richards' equation, preserving bounds of the low-order scheme.
\item Accuracy and convergence rates of the model are are verified for one- and two-dimensional test cases.
\end{keypoints}

\begin{abstract}
Simulating infiltration in porous media using Richards' equation remains computationally challenging due to its parabolic structure and nonlinear coefficients. While a wide range of numerical methods for differential equations have been applied over the past several decades, basic higher-order numerical methods often fail to preserve physical bounds on water pressure and saturation, leading to spurious oscillations and poor iterative solver convergence. Instead, low-order, bound-preserving methods have been preferred. The combination of mass lumping and relative permeability upwinding preserves bounds but degrades accuracy to first order in space. Flux-corrected transport is a high-resolution numerical technique designed for combining the bound-preserving property of low-order schemes with the accuracy of high-order methods, by blending the two methods through limited anti-diffusive fluxes. In this work, we extend flux-corrected transport schemes to the nonlinear, degenerate parabolic structure of Richards' equation, verify attainment of second-order convergence on unstructured meshes, and demonstrate applications to stormwater management infrastructure.

\end{abstract}
\section*{Plain Language Summary}
This study presents a new computational method for solving Richards' equation, a model of variably-saturated flow in soil that can simulate processes including infiltration, seepage, and water table fluctuations. The new method approximates the solution with linear functions on triangles and tetrahedra and combines two important characteristics for this representation: the optimal error convergence rate, which is quadratic, and preservation of local and global bounds on the water pressure and saturation, which ensures that numerical approximations to these quantities are physically realistic. 

\section{Introduction}
Infiltration of water into soils can be modeled by combining mass conservation with an extension of Darcy's law to variably-saturated porous media resulting in the Richards' Equation (RE) \cite{richardsCapillaryConductionLiquids1931}. Neglecting internal mass sources, we can write RE as
\begin{linenomath*}
\begin{align}
    \frac{\partial \left(\omega \rho S\right)}{\partial t} + \nabla \cdot \rho \mathbf{q} &= 0,\\
    \mathbf{q} &= - \frac{ \mathbf{\bar{k}}_i k_r}{\mu_w}\left(\nabla p - \rho \mathbf{g}\right), \label{eq:darcy}
\end{align}
\end{linenomath*}
where $\omega$ is porosity, $\rho$ is the density of water, $S$ is the water saturation, $\mathbf{q}$ is the Darcy velocity, $\mathbf{\bar{k}}_i$ is the intrinsic permeability tensor of the medium, $k_r$ is the relative permeability for the current degree of water saturation, $\mu_w$ is the dynamic viscosity of water, and $\mathbf{g}$ is the gravitational acceleration with vector magnitude $g=\| \mathbf{g}\|$. Introducing a coordinate system where the $z$ coordinate is ``positive up", meaning the unit vector in the $z$-direction, $e_z$, is $-\mathbf{g}/g$, we can introduce the gravitational potential relative to the origin, $\rho g z$ and rewrite the gravitational force term above as $-\rho g \nabla z$. 

For constant density, the introduction of the hydraulic head $\psi= p/\rho g$ and total head $\phi = \psi + z$, results in the common form of Darcy's law
\begin{linenomath*}
\begin{equation}
\mathbf{q} = - \mathbf{\bar{K}}_s k_r\nabla \phi,    
\end{equation}    
\end{linenomath*}

where $\bar{\mathbf{K}}_s = \frac{ \mathbf{\bar{k}}_i \rho g}{\mu_w}$ is the saturated hydraulic conductivity tensor (or $\bar{\mathbf{K}} = \bar{\mathbf{K}}_s k_r$ the variably-saturated hydraulic conductivity tensor). For reference, we note that a common form used in petroleum engineering and more generally multiphase flow modeling is
\begin{linenomath*}
\begin{equation}
\mathbf{q} = -  \mathbf{\bar{k}}_i \lambda \left(\nabla p - \rho \mathbf{g}\right)  
\end{equation}   
\end{linenomath*}
where $\lambda = k_r/\mu$ is the water mobility.

To close the model equations, additional pressure-saturation-relative permeability relations are required, typically resulting in closed algebraic expressions for $S$ as a function of $\psi$, and $k_r$ as a function of $S$. We may write $S(\psi)$ and $k_r(S)$ (or the composite function $k_r(\psi)$) to emphasize the choice of primary unknown. We note also that $\phi$  is a function of $\psi$, so we may likewise write $\phi(\psi)$. Introducing the volume fraction, $\theta=\omega S$, we can then write RE as
\begin{linenomath*}
\begin{equation}
\frac{\partial \theta}{\partial t} + \nabla \cdot \left( -\bar{\mathbf{K}}_s k_r \nabla \phi \right) = 0.    
\end{equation}
\end{linenomath*}
Three common formulations of RE are the $\psi$-based form, $\theta$-based form, and mixed form, which can be written as:
$\psi$-based
\begin{linenomath*}
\begin{align}
    C(\psi) \frac{\partial \psi}{\partial t}-\nabla \cdot 
    \left[\bar{\mathbf{K}} (\psi)\nabla \phi(\psi)\right]=0,
\label{eqn:h-based}
\end{align}
\end{linenomath*}
$\theta$-based
\begin{linenomath*}
\begin{align}
    \frac{\partial \theta}{\partial t}-\nabla \cdot \left[D(\theta) \nabla \theta+\bar{\mathbf{K}}(\theta) \nabla z\right]=0,
\label{eqn:theta-based}
\end{align}
\end{linenomath*}
Mixed 
\begin{linenomath*}
\begin{align}
\frac{\partial \theta}{\partial t}-\nabla \cdot \left[\bar{\mathbf{K}}(\theta)\nabla \phi (\psi) \right]=0.
\end{align}
\label{eqn:mixed-based}
\end{linenomath*}

These three distinct forms of RE help illuminate its structure. In the $\psi$-based form, the term $C(\psi) = \frac{\partial \theta}{\partial \psi}$ in equation~\eqref{eqn:h-based} is known as the specific moisture capacity and RE is a quasilinear heat equation for $\psi$ or $\phi$, suggesting the need for implicit time discretizations to avoid unreasonable time-step restrictions. Conversely, in the $\theta$-based form, the term $D(\theta) = \frac{K(\psi)}{C(\theta)}$ in equation~\eqref{eqn:theta-based} represents the unsaturated diffusivity and RE displays the structure of a degenerate advection-diffusion equation for  $\theta$, explaining the self-sharpening wetting fronts that arise and persist during ponded infiltration, and suggesting characteristics of trasport equations. The $\psi$-based form does  not conserve water mass when discretized directly \cite{celiaGeneralMassconservativeNumerical1990}, while the $\theta$-based form degenerates under the commonly occurring water-saturated state, making the \textit{mixed form} the basis for general-purpose numerical models. The choice of formulation can lead to markedly different numerical behaviors in unsaturated media with coarse discretizations (e.g. incorrect front speeds due to mass-conservation errors in the $\psi$-based form). The \textit{mixed form} of RE offers a balanced approach, combining accurate mass conservation with improved numerical performance and comparable computational cost.

Over the past few decades, much attention has been given to the development of fast and accurate numerical solvers of RE. We mention here several notable advances leading to the current approach. As mentioned above, a robust conservative numerical model was developed by combining the mixed form with mass lumping of the Galerkin finite element method \cite{celiaGeneralMassconservativeNumerical1990}. A closely related approach, which added edge-based upwinding of the relative permeability, yields a provably bound-preserving method that is first-order accurate on unstructured meshes, but the mesh must satisfy certain quality constraints to achieve these properties \cite{forsythMonotonicityConsiderationsSaturated1997}. Mass-conservative formulations of multistep methods in time achieved higher-order convergence in time but were not provably free of non-physical oscillations \cite{keesHigherOrderTime2002}. Stabilized finite element schemes for unstructured meshes were developed based on the variational multiscale approach, with the local mass conservation property recovered by post-processing, but the method required nonlinear shock-capturing to reduce non-physical oscillations, reducing the method to first-order in space \cite{keesLocallyConservativeStabilized2008}. More recently, nonlinear stability and first-order convergence for a finite element scheme similiar to the early upwinded, mass-lumped scheme or \citeA{forsythMonotonicityConsiderationsSaturated1997} was proved for RE \cite{aithammououlhajNumericalAnalysisNonlinearly2018}, and it forms a key building block of this work.

Limitations of RE as a model of the physical processes are also apparent. It can lose accuracy under saturated and ponding conditions due to neglect of air-phase effects   \cite{zhangApplicabilityRichardsEquation2026}. While there is likely a need for modeling improvements that enhance physical realism and improve predictive capability in complex environments \cite{wangAdvancingPhysicalRealism2026}, the aspects we focus on here are connected with sharp wetting fronts and the transition from variably-saturated to fully-saturated conditions, which are genuine aspects of the physical infiltration process and therefore likely challenges to enhanced models too.

Like RE, many scalar linear and nonlinear parabolic initial-boundary value problems possess maximum principles, by which we mean generally that solution extreme values are controlled by the boundary data in time and space. Many numerical methods nevertheless lack corresponding discrete bound-preservation properties for their fully discrete representation of initial-boundary value problems and hence allow the formation of new extrema, which are often recognized as non-physical oscillations and which significantly degrade the error and iterative solver performance.  In particular, standard finite element formulations often fail to enforce the discrete local maximum principle for solutions with sharp fronts \cite{kuzminAlgebraicEntropyFixes2020}. By analyzing the algebraic structure of fully-discrete linear advection equations, it can be shown that violation of the discrete maximum principle may occur when there are positive off-diagonal entries of the mass matrix \textbf{M} and negative off-diagonal entries of the discrete advection operators \cite{kuzminFluxCorrectedTransportPrinciples2012}. The same issue arises when simulating infiltration into initially very dry soils, where spurious oscillations appear in both the pressure head ($\psi$) and moisture content ($\theta$) fields, independently of whether the method conserves mass. These numerical oscillations manifest as fluctuations in the spatial distribution of pressure head around the true solution \cite{belfortEquivalentHydraulicConductivity2013, barua_bound_2023}.  

The occurrence of such oscillations, resulting from both spatial and temporal discretizations, has long been a concern in finite element simulations of unsaturated flow and numerical methods for partial differential equations more generally. When the discretized system matrix does not satisfy the so-called M-matrix condition, the numerical solution may violoate the local bounds dictated by the maximum principle~\cite{forsythMonotonicityConsiderationsSaturated1997, belfort2005comparison, younes2006new}. To enforce this condition strictly in RE, modifications to the standard Galerkin method are employed that eliminate non-negative off-diagonal entries from the discrete operator. One commonly employed modification is mass lumping, which diagonalizes the mass matrix through row summation.  However, mass lumping is not sufficient to guarantee monotonicity for RE. The first reported approach to guarantee monotonicity was to use upwinding along mesh edges with relative permeabilities defined at nodes, together with mass lumping, to obtain a provably monotone scheme under certain mesh quality assumptions \cite{forsythMonotonicityConsiderationsSaturated1997}. Despite the long history of work, rigorous proofs of nonlinear stability and positivity preservation for the lumped upwind finite element discretization of RE were not available until relatively recently \cite{aithammououlhajNumericalAnalysisNonlinearly2018}.

This work aims at extending these well-known and rigorous first-order approaches to bound-preserving finite element discretizations of RE to higher order methods. The flux-corrected transport (FCT) methodology introduced by \citeA{borisFluxcorrectedTransportSHASTA1973} in the early 1970s enables the removal of excessive artificial diffusion from a bound-preserving low-order predictor in a local extremum diminishing manner. Zalesak's multidimensional generalization of the FCT paradigm extends its applicability to unstructured meshes and finite elements \cite{zalesak1979fully}. It is employed in many representatives of modern Algebraic Flux Correction (AFC) schemes for hyperbolic and convection-dominated conservation laws. Comprehensive reviews of FCT/AFC methods and their theoretical foundations can be found in \cite{kuzminFluxCorrectedTransportPrinciples2012, kuzmin2024property}. Although flux limiters can, in principle, be incorporated into high-order schemes for the RE, the nonlinear dependence of the moisture content on the pressure head in the mixed formulation precludes the direct application of existing limiting techniques.

In this work, we present a new scheme that uses the reliable mixed form of RE discretized by the standard linear Galerkin FEM with algebraic stabilization based on mass lumping and upwinding as the basic building block. Then, to achieve optimal accuracy while maintaining the bound-preserving property of the first-order building block, we develop a dedicated FCT algorithm for the mixed form of RE. Specifically, the low-order component of our scheme is the nonlinear upwind method analyzed in \cite{aithammououlhajNumericalAnalysisNonlinearly2018}. The difference between the high- and low-order approximations of the moisture content is decomposed into antidiffusive fluxes, which are constrained using Zalesak’s multidimensional limiter. The main novelty of our approach lies in combining FCT with implicit time discretizations for a nonlinear problem that may exhibit locally elliptic behavior. Another distinctive feature is the iterative recovery of the pressure head from the flux-corrected moisture content. The preservation of local bounds for the latter is ensured by the limiting procedure. Importantly, the flux correction step remains fully explicit, although nonlinear systems must be solved iteratively in other stages of the algorithm, and the method is comparable in cost to the first-order, bound-preserving scheme. Finally, the practical application of this modeling approach is demonstrated through the simulation of infiltration patterns in stormwater infiltration systems.

\section{Fully Discrete Galerkin Method for RE}

In this work, we evolve the pressure head $\psi$ using the mixed (conservative)
form 
\begin{linenomath*}
  \begin{align}
    \frac{\partial \theta(\psi)}{\partial t} - \nabla \cdot \left(\bar{\mathbf{K}}(\psi)(\nabla \psi + \nabla z) \right) = 0
   \label{eqn:MixedForm}
\end{align}
\end{linenomath*}
of the RE with $\bar{\mathbf{K}}(\psi)=
\bar{\mathbf{K}}_s(\mathbf x)k_r(\psi)$. Both the porous medium and the fluid are assumed to be incompressible. The hydraulic conductivity, moisture content, total head, and other parameters are nonlinear functions of $\psi$. The effective saturation is defined by 
\begin{linenomath*}
    \begin{equation}
        S_e(\psi)=\frac{\theta(\psi)-\theta_r}{\theta_s-\theta_r},        
    \end{equation}
\end{linenomath*}
where $\theta_r$ and $\theta_s$ denote the residual and saturated moisture contents, respectively. Solving the above equation for $\theta(\psi)$ yields
$\theta(\psi) = \theta_r+ S_e(\psi)(\theta_s-\theta_r)$.
Standard constitutive relations linking the pressure head to soil properties are summarized in Table~\ref{tab:psk}. Among these, the van Genuchten model is the most widely used closure for the RE, largely due to its ability to produce smooth constitutive profiles.

\begin{table}[!ht]
\caption{\label{tab:psk} Pressure-saturation-conductivity relations.}
\scalebox{0.8}{
\begin{tabular}{m{2cm} m{7cm} m{5cm}}
\hline
\textbf{Model} & \textbf{Formulation} & \textbf{Parameters} \\ 
\hline

Brooks--Corey Model 
& 
\begin{tabular}{llp{.15\textwidth}}
$S_e(\psi)=
\begin{cases}
\alpha |\psi|^{-n}, & \psi < -\dfrac{1}{\alpha} \\
1, & \psi \ge -\dfrac{1}{\alpha}
\end{cases}$ \\[0.2cm]
$k_r(\psi)= S_e^{L + 2 + \frac{2}{n}}$
\end{tabular}
&
\begin{tabular}{llp{.15\textwidth}}
$\alpha$: inverse air-entry value \\
$n$: pore-size index \\
$L$: pore-connectivity parameter \\
$L = 0.5$--$1$
\end{tabular}
\\ \hline

Haverkamp et al. Model 
&
\begin{tabular}{llp{.15\textwidth}}
$S_e(\psi)=
\begin{cases}
\dfrac{B}{B + |\psi|^{\beta}}, & \psi < 0 \\
1, & \psi \ge 0
\end{cases}$ \\[0.2cm]
$k_r(\psi)= \dfrac{A}{A + |\psi|^{\alpha}}$
\end{tabular}
&
\begin{tabular}{llp{.15\textwidth}}
$A, B$: scaling parameters \\
$\alpha, \beta$: empirical exponents
\end{tabular}
\\ \hline

van Genuchten Model
&
\begin{tabular}{c}
$
S_e(\psi)=
\begin{cases}
\left[1+\alpha|\psi|^n\right]^{-m}, & \psi < 0 \\
1, & \psi \ge 0
\end{cases}
$ \\[0.2cm]
$
k_r(\psi)= S_e^{L}
\left[1-(1-S_e^{1/m})^m\right]^2
$
\end{tabular}
&
\begin{tabular}{llp{.15\textwidth}}
$n, m$: pore-size parameters \\
$L$: pore-connectivity parameter \\
$L = 0.1$--$0.5$
\end{tabular}
\\ \hline

Modified van Genuchten Model
&
\begin{tabular}{llp{.05\textwidth}}
$
S_e(\psi)=
\begin{cases}
\left[S_E^*\left(1+\alpha|\psi|^n\right)^m\right]^{-1},
& \psi < -\psi_e \\
1, & \psi \ge -\psi_e
\end{cases}
$ \\[0.2cm]
$
k_r(\psi)= S_e^{L}
\left[
\dfrac{1-(1-(S_E^* S_e)^{1/m})^m}
{1-(1-S_E^{1/m})^m}
\right]^2
$
\end{tabular}
&
\begin{tabular}{llp{.05\textwidth}}
$S_E^*=\left[1+(\alpha \psi_e)^n\right]^m$ \\
$\psi_e$: equivalent air-entry value
\end{tabular}
\\ \hline

Gardner Model
&
\begin{tabular}{llp{.15\textwidth}}
$
S_e(\psi)=
\begin{cases}
\exp(\alpha \psi), & \psi < 0 \\
1, & \psi \ge 0
\end{cases}
$ \\[0.2cm]
$k_r(\psi)= S_e$
\end{tabular}
&
\begin{tabular}{llp{.15\textwidth}}
$\alpha$: inverse capillary length \\
controls exponential decay
\end{tabular}
\\ \hline

\end{tabular}}
\end{table}

\begin{figure}[ht]
\centering
\includegraphics[width=0.7\textwidth]{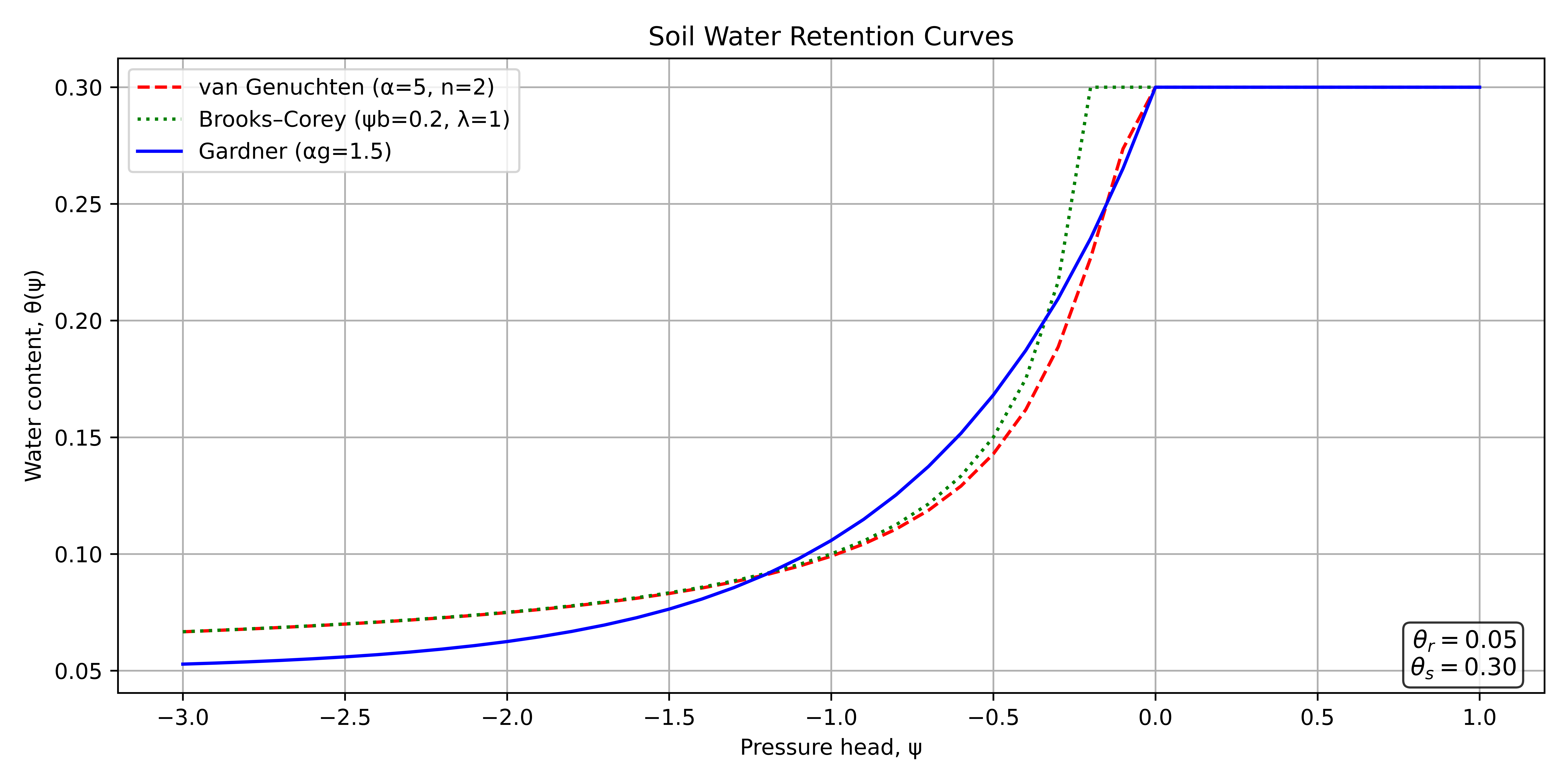}
\caption{Pressure head vs moisture content. }
\label{fig:psk_comparison}
\end{figure}

Figure ~\ref{fig:psk_comparison} shows the relationship between pressure head $\psi$ and moisture content $\theta$ for different constitutive relations  in Table~\ref{tab:psk}. The plotting is done for $\theta_r = 0.05$ and $\theta_s = 0.30$. Different relations show different smoothness near saturation.

Restricting our attention to a bounded domain $\Omega$, we complete the problem formulation by prescribing appropriate initial and boundary conditions:
\begin{align}\label{richBC}
  \psi|_{t=0} = \psi_0,\qquad
  \psi|_{\partial\Omega_D} = \psi_D,\qquad 
  -\bar{\mathbf{K}}(\psi)\nabla \phi \cdot \mathbf{n}
  |_{\partial\Omega_N}= g_N,
\end{align}
where
$\mathbf{n}$ is the unit outward normal to the boundary
$\partial\Omega=\partial\Omega_D\cup\partial\Omega_N$ and
\begin{linenomath*}
    \begin{equation}
        \phi=\psi+z.        
    \end{equation}
\end{linenomath*}

To prepare for a finite element discretization of our initial-boundary-value problem, we multiply \eqref{eqn:MixedForm} by a test function  $w\in H^1(\Omega)$ satisfying the homogeneous Dirichlet condition $w|_{\partial\Omega_D}=0$. The resulting weighted residual is integrated over $\Omega$, and integration by parts is applied to the diffusive term. The boundary term 
\begin{linenomath*}
    \begin{equation}
        \int_{\partial \Omega} w \left[\bar{\mathbf{K}}(\psi)\nabla\phi \cdot \mathbf{n} \right] \, dS,        
    \end{equation}
\end{linenomath*}

is then decomposed into integrals over the Dirichlet and Neumann portions of $\partial\Omega$. The integral over $\partial\Omega_D$ vanishes due to the choice of the test space. The integral over $\partial\Omega_N$ incorporates the prescribed Neumann data into the weak formulation
\begin{linenomath*}
\begin{equation}\label{weakform}
\int_{\Omega} w \frac{\partial \theta(\psi)}{\partial t} \, dV
+ \int_{\Omega} \bar{\mathbf{K}}(\psi) \nabla \phi \cdot \nabla w \, dV
= -\int_{\partial \Omega_N} w g_N \, dS,
\end{equation}
\end{linenomath*}
which is used in the standard Galerkin finite element method (FEM) for
solving \eqref{eqn:MixedForm}.

\subsection{Space Discretization}

To discretize the Galerkin weak form \eqref{weakform} of the RE in space, we use a conforming finite element partition of $\Omega$ with nodal points $\{x_j\}$. The numerical approximations
\begin{linenomath*}
    \begin{equation}
        \psi_h = \sum_j \psi_j v_j,\qquad
        \phi_h = \sum_j \phi_j v_j,\qquad w_h = \sum_i w_i v_i
    \end{equation}
\end{linenomath*}
are defined using continuous, piecewise-linear basis functions $v_i$ that
possess the nodal interpolation property
$v_i(x_j)=\delta_{ij}$. Following the \emph{group finite element}
approach analyzed in \citeA{barrenecheaAnalysisGroupFinite2017}, we
approximate the nonlinear function
$\theta(\psi_h)$ by its nodal interpolant,
$\theta_h = \sum_j \theta(\psi_j) v_j$. This treatment is
equivalent to using
inexact nodal quadrature for the integral involving the time derivative
of $\theta(\psi_h)$.

The semi-discrete weak form using $w=v_i$ and the above approximations reads
\begin{linenomath*}
\begin{equation}
    \sum_j \left[\int_{\Omega} v_j v_i \, dV \right] \frac{d \theta_j}{dt}
    + \sum_j \left[\int_{\Omega} \bar{\mathbf{K}}(\psi_h) \nabla v_j \cdot \nabla v_i \, dV \right] \phi_j \notag
    + \sum_j \left[\int_{\partial\Omega_N} g_N v_i \, dS \right] = 0.    
\end{equation}
\end{linenomath*}
The nonlinear system of evolution equations for the time-dependent nodal
values $\psi_i=\phi_i-z_i=\psi(\theta_i)$ of the pressure head can be written in the matrix form
\begin{linenomath*}
\begin{equation}\label{nonlinsys}
\mathbf{M} \frac{d \boldsymbol{\theta}(\boldsymbol{\psi})}{dt} + \mathbf{A}(\boldsymbol{\psi}) \boldsymbol{\phi}(\boldsymbol{\psi}) + \mathbf{b}^{(\text{bc})} = 0,
\end{equation}    
\end{linenomath*}
where \( \mathbf{M}=(m_{ij}) \) is the consistent mass matrix, \( \mathbf{A}=(a_{ij} )\) is the conductivity (stiffness) matrix, \( \mathbf{b}^{(\text{bc})} =(b^{bc}_i  )\) is the contribution of the Neumann boundary term, and $\boldsymbol{\psi}=(\psi_i)$ is the vector of degrees of freedom. The coefficients of \eqref{nonlinsys} are given by
\begin{align*}
    m_{ij} = \int_{\Omega} v_j v_i \, dV,\qquad
    a_{ij} (\boldsymbol{\psi})= \int_{\Omega} \bar{\mathbf{K}}(\psi_h) \nabla v_j \cdot \nabla v_i \, dV,\qquad
    {b}^{bc}_i =  \int_{\partial\Omega_N} g_N v_i \, dS
\end{align*}
and can be computed using the standard element-by-element assembly procedure.

\subsection{Temporal Discretization}

To avoid stability problems due to the elliptic behavior of the
RE in saturated regions, we discretize \eqref{nonlinsys} in time using
an implicit version of the  \( \Theta \)-scheme
\begin{equation}\label{thetaNL}
\mathbf{M} \frac{\boldsymbol{\theta}(\boldsymbol{\psi}^{n+1}) - \boldsymbol{\theta}(\boldsymbol{\psi}^n)}{\Delta t}
+ \Theta \mathcal{L}(\boldsymbol{\psi}^{n+1}) + (1 - \Theta) \mathcal{L}(\boldsymbol{\psi}^n) = 0,
\end{equation}
where $\mathcal{L}(\boldsymbol{\psi}) = \mathbf{A}(\boldsymbol{\psi}) \boldsymbol{\phi}(\boldsymbol{\psi}) + \mathbf{b}^{(\text{bc})}$ denotes the
nonlinear operator of the spatial semi-discretization. The parameter \( \Theta \in [0, 1] \) controls the degree of implicitness. 


\subsection{Newton Iteration}

The nonlinear system \eqref{thetaNL}
arising from an implicit Galerkin discretization of the RE requires an iterative solution because both the moisture content \( \theta \) and the hydraulic conductivity \( \bar{\mathbf{K}} \) depend nonlinearly on the pressure head \( \boldsymbol{\psi} \). The Newton method updates the current iterate \( \boldsymbol{\psi}^{n+1}_{k} \) via the linearized relation
\begin{linenomath*}
\begin{equation}
    \boldsymbol{\psi}^{n+1}_{k+1} = \boldsymbol{\psi}^{n+1}_{k} - \mathbf{J}^{-1}(\boldsymbol{\psi}^{n+1}_{k}) \mathcal{R}(\boldsymbol{\psi}^{n+1}_{k}).
\end{equation}
\end{linenomath*}
The Jacobian \( \mathbf{J}(\boldsymbol{\psi})\) is obtained by differentiating the nonlinear residual
\begin{linenomath*}
\begin{equation}\label{residH}
    \mathcal{R}(\boldsymbol{\psi}) =
    \mathbf{M} \frac{\theta(\boldsymbol{\psi}) - \theta(\boldsymbol{\psi}^n)}{\Delta t}
    +\Theta \mathcal{L}(\boldsymbol{\psi}) + (1 - \Theta) \mathcal{L}(\boldsymbol{\psi}^n)
\end{equation}
\end{linenomath*}
with respect to \( \boldsymbol{\psi} \), giving
\begin{linenomath*}
\begin{equation}
    \mathbf{J}(\boldsymbol{\psi}) =
      \frac{1}{\Delta t} \mathbf{M}
      \frac{d \theta(\boldsymbol{\psi})}{d \boldsymbol{\psi}}
    +\Theta \frac{\partial \mathcal{L}(\boldsymbol{\psi})}{\partial \boldsymbol{\psi}}.
\end{equation}
\end{linenomath*}
The iteration proceeds until a specified tolerance is reached for a norm
 of $\mathcal{R}(\boldsymbol{\psi}^{n+1}_{k})$.

 \section{First-Order, Boundary-Preserving Method}

 The standard Galerkin discretization of the RE requires special stabilization due to the highly nonlinear and spatially varying nature of the hydraulic conductivity \( \bar{\mathbf{K}}(\psi)\). A~nonlinearly stable and positive low-order approximation to \eqref{thetaNL} can be constructed using row-sum mass lumping for the time derivative term and nodal upwinding for the elliptic part. Following \cite{forsythMonotonicityConsiderationsSaturated1997} and \cite{aithammououlhajNumericalAnalysisNonlinearly2018}, we modify the discrete operators of
the Galerkin space discretization \eqref{nonlinsys} as follows:
\begin{linenomath*}
    \begin{equation}
        m_{ij}\approx m_i\delta_{ij},\qquad
        a_{ij}\approx\tilde a_{ij}=\tilde k_{r,ij}
        \int_{\Omega} \mathbf{\bar{K}_s}\nabla v_j \cdot \nabla v_i \, dV=k_{r,ij} a_{s,ij},        
    \end{equation}
\end{linenomath*}
where $m_i=\sum_jm_{ij}$ is a diagonal entry of the lumped mass matrix $\tilde{\mathbf{M}}=(\tilde m_{ij})$ and
\begin{linenomath*}
    \begin{equation}
        \tilde k_{r,ij} =
    \begin{cases}
    k_r(\psi_i) & \text{if } \, a_{s,ij}(\phi_j - \phi_i) \geq 0,\\
    k_r(\psi_j) & \text{otherwise}
    \end{cases}
    \end{equation}
\end{linenomath*}
is the upwind hydraulic conductivity depending on the entries $k_{ij}$ of $\bar{\mathbf{K}}_s=(k_{ij})$. The entries $\tilde a_{ij}$ constitute the upwind
conductivity matrix $\tilde{\mathbf{A}}=(\tilde a_{ij})$. In the case of strong anisotropy in $\mathbf{\bar{K}}_s$ or poorly shaped (anisotropic) elements, the above method can nonetheless fail to preserve local bounds, as already recognized \cite{forsythMonotonicityConsiderationsSaturated1997}. In such cases, we can make the further modification, for $i\neq j$,
\begin{linenomath*}
    \begin{equation}
        \tilde{a}_{ij} =-\tilde{k}_{rij} \max(0,-a_{s,ij}),        
    \end{equation}
\end{linenomath*}

and finally
\begin{linenomath*}
    \begin{equation}
        \tilde{a}_{ii} = -\sum_{j \neq i} \tilde{a}_{ij},    
    \end{equation}
\end{linenomath*}
enforcing the M-matrix condition at the cost of local inconsistency \cite{kuzminConstrainedFiniteElement2009}.  Note, that \citeA{aithammououlhajNumericalAnalysisFinite2018} show that such a correction is not strictly necessary for convergence and {\em global} bound preservation but mention that their tests show that strong anisotropy ratios produce excessive numerical diffusion and slow convergence.
The nonlinear residual of the fully discrete low-order scheme is given by
\begin{linenomath*}
\begin{equation}\label{residL}
\tilde{\mathcal{R}}(\boldsymbol{\psi}) =
\mathcal{R}(\boldsymbol{\psi}) -\mathbf f(\boldsymbol{\psi}),
\end{equation}    
\end{linenomath*}
where 
\begin{linenomath*}
\begin{align}\nonumber
  \mathbf f(\boldsymbol{\psi})=&\frac{1}{\Delta t}
  (\tilde{\mathbf{M}}-\mathbf{M})
  (\boldsymbol{\psi}-\boldsymbol{\psi}^n)
  +\Theta[\tilde{\mathbf{A}}(\boldsymbol{\psi})-\mathbf{A}(\boldsymbol{\psi})]
  \boldsymbol{\phi}(\boldsymbol{\psi}) \\
  &+(1-\Theta)[\tilde{\mathbf{A}}(\boldsymbol{\psi}^n)-\mathbf{A}(\boldsymbol{\psi}^n)]\boldsymbol{\phi}(\boldsymbol{\psi}^n).\label{fcorr}
\end{align}
\end{linenomath*}
The removal of $\mathbf f(\boldsymbol{\psi})$
from the residual $\mathcal{R}(\boldsymbol{\psi})$ defined by
\eqref{residH} introduces significant amounts of artificial diffusion, which stabilizes the Galerkin discretization and prevents spurious oscillations in regions with steep solution gradients. Moreover, the analysis performed in \cite{forsythMonotonicityConsiderationsSaturated1997,aithammououlhajNumericalAnalysisNonlinearly2018} establishes the validity of discrete maximum principles for both $\theta$ and $\phi$ under appropriate assumptions.

\section{Second-Order Flux-Corrected Transport Step}
Let $\boldsymbol{\psi}^H$ and $\boldsymbol{\psi}^L$ denote solutions of the high-order system $\mathcal{R}(\boldsymbol{\psi}^H)=0$ and of the low-order system $\tilde{\mathcal{R}}(\boldsymbol{\psi}^L)=0$, respectively. 
In view of \eqref{residL} and  \eqref{fcorr}, we have
\begin{linenomath*}
\begin{equation}
    \tilde{\mathbf{M}}\boldsymbol{\theta}(\boldsymbol{\psi}^H)
=\tilde{\mathbf{M}}\boldsymbol{\theta}(\boldsymbol{\psi}^L)
+\Delta t\mathbf{f}(\boldsymbol{\psi}^H,\boldsymbol{\psi}^L,\boldsymbol{\psi}^n).
\end{equation}    
\end{linenomath*}

The matrices $\tilde{\mathbf{A}}$, $\mathbf{A}$ and
$\tilde{\mathbf{M}}-\mathbf{M}$ that appear in \eqref{fcorr}
are symmetric with zero row and column sums. Therefore, they represent discrete diffusion (alias \emph{graph Laplacian}) operators and the corresponding matrix-vector products can be decomposed into numerical fluxes associated
with edges of the sparsity graph \cite{kuzminFluxCorrectionTools2002,kuzminFluxCorrectedTransportPrinciples2012,kuzmin2024property}. A flux decomposition of $\mathbf f=(f_i)$ is given by

\begin{linenomath*}
\begin{align}
f_i=\sum_{j\ne i}f_{ij},\qquad f_{ji}=-f_{ij}.    
\end{align}
\end{linenomath*}

A classical FCT algorithm multiplies the flux $f_{ij}$ by an adaptively chosen correction factor $\alpha_{ij}\in [0,1]$. The flux $f_{ji}$ is limited using $\alpha_{ji}=\alpha_{ij}$. The purpose of this flux correction is to
enforce local discrete maximum principles of the form
\begin{linenomath*}
\begin{equation}\label{FCTconstraints}
  \theta_i^{\min}\le \theta_i^{n+1}=\theta_i^L+\frac{\Delta t}{m_i}\sum_{j\ne i}
  \alpha_{ij}f_{ij}\le \theta_i^{\max}.
\end{equation}    
\end{linenomath*}
The bounds of these  constraints are constructed using
$\theta_j^n=\theta(\psi_j^n)$ and/or  $\theta_j^L=\theta(\psi_j^L)$. By default,
$\theta_i^{\min}=\min_{j\in\mathcal N_i}\theta_j^L,\ \theta_i^{\max}=\max_{j\in\mathcal N_i}\theta_j^L$,
where $\mathcal N_i=\{j\,:\, m_{ij}>0\}$.

If the Crank--Nicolson time discretization corresponding to the choice
$\Theta=\frac12$ is used not only for the
high-order system but also for the low-order one, then 
\begin{linenomath*}
\begin{align}
  f_{ij}=&m_{ij}\left(\frac{\theta_i^{H}-\theta_i^{n}}{\Delta t}-
\frac{\theta_j^{H}-\theta_j^{n}}{\Delta t}\right)
-\frac12 a_{ij}(\boldsymbol{\psi}^H)(\phi(\psi_j^H)-\phi(\psi_i^H)) \nonumber\\
&+\frac12\tilde a_{ij}(\boldsymbol{\psi}^L)(\phi(\psi_j^L)-\phi(\psi_i^L))
-\frac12(a_{ij}(\boldsymbol{\psi}^n)-
\tilde a_{ij}(\boldsymbol{\psi}^n))(\phi(\psi_j^n)-\phi(\psi_i^n)).
\end{align}
\end{linenomath*}
However, the semi-implicit treatment of the low-order system introduces severe time
step restrictions even in the case of linear transport problems \cite{kuzminFluxCorrectionTools2002,kuzminFluxCorrectedTransportPrinciples2012,kuzmin2024property}. The nonlinear stability and positivity of the low-order upwind method for the RE were shown in \cite{forsythMonotonicityConsiderationsSaturated1997,aithammououlhajNumericalAnalysisNonlinearly2018} only for the fully implicit (backward Euler, $\Theta=1$) version. 
In the present work, both the high-order and low-order transport operators are discretized in time by backward Euler, i.e. $\Theta=1$. First, the low-order solution $\boldsymbol{\psi}^L$ is computed from the monotone nonlinear system. An explicit FCT correction is applied to the low-order update in terms of the conserved variable $\theta$. For the backward Euler / backward Euler combination used here, the edge-based raw antidiffusive fluxes take the form
\begin{linenomath*}
\begin{equation}
    f_{ij}
=
m_{ij}\left(
\frac{\theta_i^H-\theta_i^n}{\Delta t}
-
\frac{\theta_j^H-\theta_j^n}{\Delta t}
\right)
-
a_{ij}(\boldsymbol\psi^H)(\phi_j^H-\phi_i^H)
+
\tilde a_{ij}(\boldsymbol\psi^L)
\end{equation}    
\end{linenomath*}

The multidimensional limiting strategy of \citeA{zalesak1979fully}  enforces the FCT constraints \eqref{FCTconstraints} using correction factors $\alpha_{ij}$
such that the sum of positive limited fluxes cannot create an overshoot, while the 
sum of limited fluxes cannot create an undershoot. Maximization of $\alpha_{ij}\in[0,1]$
under these constraints leads to the following algorithm.

\begin{algorithm}[ht]
\caption{Zalesak's multidimensional FCT limiter.}

\begin{algorithmic}[1]
\For{each node \( i \)}
  \State Compute: \( P_i^+ = \sum_{j\ne i} \max(f_{ij}, 0) \), \quad \( P_i^- = \sum_{j\ne i} \min(f_{ij}, 0) \)
  \State 
  \(Q_i^+ = \frac{m_i}{\Delta t}(\theta_i^{\max} - \theta_i^L),
  \quad Q_i^- = \frac{m_i}{\Delta t}(\theta_i^{\min} - \theta_i^L)\)
  \State 
  \(
    R_i^+ = \min\left(1, \frac{Q_i^+}{P_i^+} \right), \quad R_i^- = \min\left(1, \frac{Q_i^-}{P_i^-} \right)
  \)
\EndFor

\For{each node pair \( (i, j) \)}
  \If{\( f_{ij} > 0 \)}
    \State \( \alpha_{ij} = \min(R_i^+, R_j^-) \)
  \ElsIf{\( f_{ij} < 0 \)}
    \State \( \alpha_{ij} = \min(R_i^-, R_j^+) \)
  \Else
    \State \( \alpha_{ij} = 1 \)
  \EndIf
\EndFor
\end{algorithmic}
\end{algorithm}

The main steps of our FCT scheme for RE \eqref{eqn:MixedForm} are as follows:
\begin{enumerate}
\item Solve the low-order nonlinear system $\tilde{\mathcal{R}}(\boldsymbol{\psi}^L)=0$
  for $\boldsymbol{\psi}^L$.
\item Compute the raw anti-diffusive fluxes $f_{ij}$ and the correction factors $\alpha_{ij}$.
\item Assemble the components $\bar f_i=\sum_{j\ne i}\alpha_{ij}f_{ij}$ of
  $\bar{\mathbf f}=(\bar f_i)$.
\item Compute the moisture content
$\boldsymbol{\theta}^{n+1}
=\boldsymbol{\theta}^L
+\Delta t \tilde{\mathbf{M}}^{-1}
\bar{\mathbf{f}}$.
\item Recover a pressure head $\boldsymbol{\psi}^{n+1}$ satisfying
  $\boldsymbol{\theta}(\boldsymbol{\psi}^{\rm n+1})=\boldsymbol{\theta}^{n+1}$.
\end{enumerate}  

It remains to discuss Step 5. In principle, in the unsaturated portion of the domain, a nodal moisture content
\( \theta_i^{n+1} \) can be converted to a pressure head \( \psi_i^{n+1} \) using the van Genuchten formula
\begin{linenomath*}
\begin{equation}
\label{eqn:VG-FCT}
\psi(\theta) = -\frac{1}{\alpha}
\left[ \left( \frac{\theta - \theta_r}{\theta_s - \theta_r} \right)^{-\frac{1}{m}} - 1 \right]^{\frac{1}{n}}. 
\end{equation}
\end{linenomath*}
The proposed FCT scheme guarantees that the flux-corrected moisture content $\theta$ satisfies
\( \theta_r < \theta < \theta_s \). This 
ensures consistency with the soil water retention model.
However, the nonlinear mapping $\psi(\theta)$ is ill-conditioned, especially in the
limits $\theta\to\theta_r$ and $\theta\to \theta_s$. As a consequence, small numerical
errors in $\theta$ may cause large jumps in $\psi$, and since the flux of the
RE depends on $\nabla\psi$, noise in $\psi$ creates noise in the flux.

In our numerical implementation, the conversion of the moisture content into a pressure head is carried out using a Newton iteration rather than the analytical formula \eqref{eqn:VG-FCT}.  The iterative approach is more stable and extends naturally to the compressible case when the density $\rho = \rho(\psi)$ is a nonlinear function of pressure. Moreover, physics-based regularization can be introduced via artificial compressibility.
The nonlinear problems to be solved are of the form $g(\psi^{n+1})=0$. The residual 
\begin{linenomath*}
\begin{equation}
g(\psi)=\rho(\psi)\,\theta(\psi)-\rho(\psi^L) \theta^{n+1}\qquad      
\end{equation}
\end{linenomath*}
depends on the density closure $\rho(\psi)$. The corresponding Jacobian is 
\begin{linenomath*}
\begin{equation}
g'(\psi)
=
\rho(\psi)\left[
\beta\,\theta(\psi)
+
\frac{d\theta}{d\psi}(\psi)
\right],
\end{equation}
\end{linenomath*}
where $\beta = \frac{d\rho}{d\psi}(\psi)/\rho$ is the water compressibility. A natural initial guess for pressure head recovery via the Newton iteration
\begin{linenomath*}
\begin{equation}
\label{eq: newton_m_2_u}
\psi_{(k+1)}^{n+1}  
=
\psi_{(k)}^{n+1} 
-
\frac{g\!\left(\psi^{n+1}_{(k)}\right)}{g'\!\left(\psi^{n+1}_{(k)}\right)},\qquad k=0,1,2,\ldots
\end{equation}
\end{linenomath*}
is the low-order predictor $\psi^{n+1}_{(0)}=\psi^L$ that defines the moisture content $\theta^L=\theta(\psi^L)$.


\section{Numerical Results}
To evaluate the performance and robustness of the proposed Flux-Corrected Transport (FCT) scheme for variably saturated flow, several numerical experiments are conducted. These include standard one-dimensional benchmark problems from the literature, as well as two- and three-dimensional simulations of engineered green infrastructure systems, such as rain gardens and bioswales. The FCT scheme is implemented within the Proteus computational framework \cite{proteus2023} and applied to the mixed-form RE to ensure mass conservation and eliminate non-physical oscillations.

Several multidimensional test problems are simulated. For the convergence test, grid refinement studies are solved by comparing the numerical solution with the analytical solution developed by \citeA{https://doi.org/10.1029/2005WR004638}.  At first, three one-dimensional test problems are conducted under varying initial conditions, such as static and hydrostatic conditions. The Van Genuchten parameters and boundary conditions are also described for each test problem to demonstrate the numerical capabilities of different numerical schemes. The Van Genuchten soil parameters used for one-dimensional test problems are summarized in Table~\ref{tab:VGparams_combined}.
 
\begin{table}[ht]
    \centering
    \caption{Van Genuchten soil parameters used in one dimensional test problems.}
    \label{tab:VGparams_combined}
    \scalebox{0.85}{
    \begin{tabular}{llcccccc}
    \textbf{Test Case} & \textbf{Soil Type} & $K_s$ (m/day) & $\theta_r$ & $\theta_s$ & $n$ & $m$ & $\alpha$ (1/m) \\ \hline
    Celia et al. (1D)~\cite{celiaGeneralMassconservativeNumerical1990} & Sand & 7.967 & 0.102 & 0.368 & 2.00 & -- & 3.35 \\
    Test Problem 2 (1D) & Sand & 7.128 & 0.045 & 0.43 & 2.68 & 0.627 & 14.5 \\
    Test Problem 3 (1D) & Sand & 0.72 & 0.045 & 0.43 & 1.592 & 0.372 & 1.389 \\ \hline
    \end{tabular}
    }
\end{table}
\subsection{Test Problem 1: Celia Test}
\begin{figure}[h]
\centering
\includegraphics[width=0.6\textwidth]{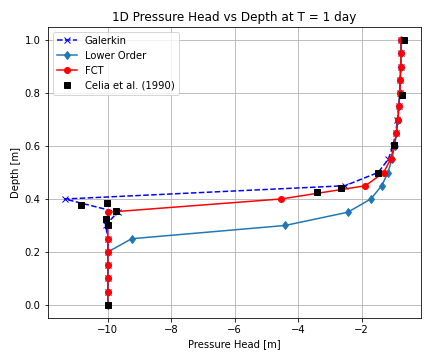}
\caption{Pressure head at 1 day (Celia et al.). }
\label{fig:result_celia}
\end{figure}

This test problem is collected from  \citeA{celiaGeneralMassconservativeNumerical1990} which is widely used for validating numerical solutions of RE. Here, the simulation considers one dimensional vertical infiltration through 100 cm sand column. Dirichlet boundary condition is applied at the top of the sand column with a pressure head of $\Psi_{top}$ = -75 cm and at the bottom $\Psi_{bottom}$ = -1000 cm. The initial condition of the pressure head is uniformly set to -1000 cm throughout the whole sand column. Here, infiltration for 60 minutes is simulated.  The simulation of infiltration of the test problem is shown in Figure \ref{fig:result_celia}. 

Different schemes are used to simulate the infiltration in the sand column. The soil properties of the sand column are provided in Table \ref{tab:VGparams_combined}. The figure shows that the standard Galerkin method creates oscillations.

\subsection{Test Problem 2: Hydrus 1D Column Comparison}
This test problem models vertical infiltration in a homogeneous sand column with a height of 20~m. The objective is to validate the numerical model against a well-defined benchmark by observing the temporal evolution of pressure head, hydraulic conductivity, and moisture content. The soil hydraulic properties used for this test case follow the Van Genuchten-Mualem formulation and are summarized in Table~\ref{tab:VGparams_combined}.

\begin{figure}[H]
\centering
a)\includegraphics[width=0.3\textwidth]{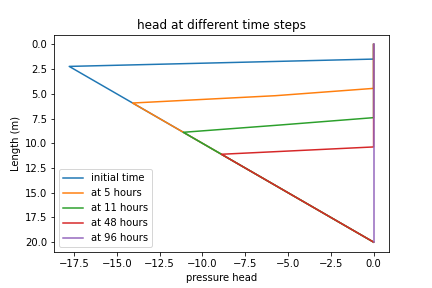} 
b)\includegraphics[width=0.3\textwidth]{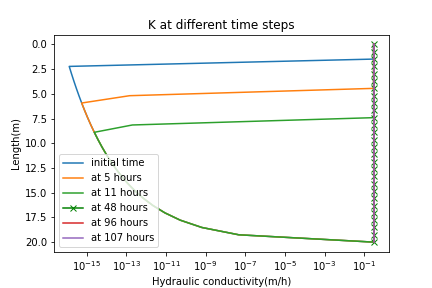} 
c)\includegraphics[width=0.3\textwidth]{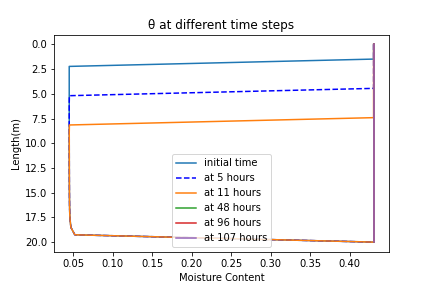}
\caption{Evolution of the pressure head (a), hydraulic conductivity (b), and
moisture content (c).}
\label{fig:test2}
\end{figure}

The initial pressure head distribution is set to hydrostatic conditions with a uniform boundary condition of $\Psi_{\text{top}} = \Psi_{\text{bottom}} = 0$. Figure~\ref{fig:test2} shows the evolution of (a) pressure head, (b) hydraulic conductivity, and (c) moisture content at different time steps. The profiles demonstrate the propagation of the wetting front and the nonlinear behavior of hydraulic properties as a function of pressure.

\begin{figure}[H]
\centering
\includegraphics[width=0.3\textwidth]{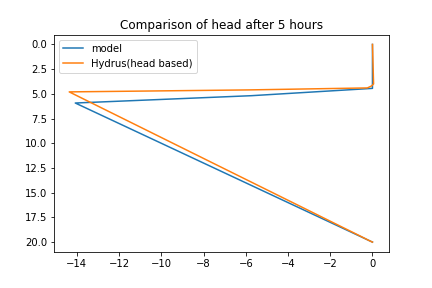} 
\includegraphics[width=0.3\textwidth]{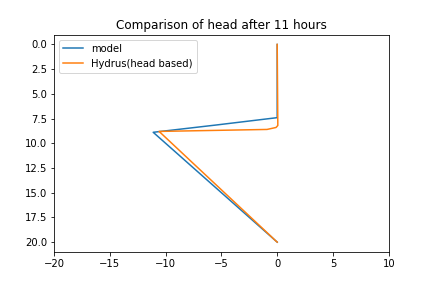}
\includegraphics[width=0.3\textwidth]{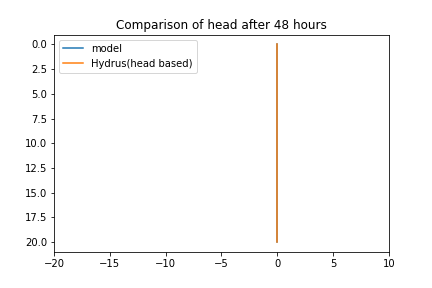}
\caption{Pressure head at 5 hours (a), 11 hours (b), and 48 hours (c).}
\label{fig:pressure2}
\end{figure}

The numerical simulation is obtained using the standard Galerkin method with mass lumping. These results were compared with the results of the HYDRUS tool. The pressure head distribution looks similar to the HYDRUS result. The comparison is shown in Figure \ref{fig:pressure2}. Here, the sand column gets saturated after 48 hours. 

\subsection{Test problem 3: Szymkiewicz 1D Sand Column}
\begin{figure}[ht]
\centering
\includegraphics[width=0.6\textwidth]{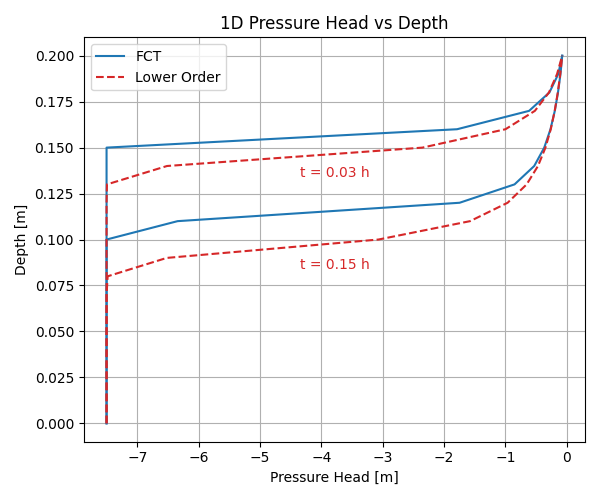}
\caption{ Comparison of different schemes for test problem 3. }
\label{fig:testc3}
\end{figure}

In this test problem, infiltration in a 20-cm-deep sand column is investigated. The initial condition assumes a uniform pressure head distribution of $\Psi(z, 0) = -750$~cm throughout the domain. The initial pressure distribution in the sand column is uniform. The initial pressure head is $\Psi$(z,0) = $\Psi_{init}$ = -750 cm. No influx and outflow at the boundary are considered for this test problem. he Dirichlet boundary conditions $\Psi_{top} = -7.5$~cm and $\Psi_{bottom} = -750$~cm are prescribed at the top and bottom of the soil column, respectively. This test problem is taken and modified from \citeA{szymkiewiczApproximationInternodalConductivities2009}. The infiltration in the homogeneous sand column is simulated for 0.15 hours. The corresponding Van Genuchten properties are provided in Table \ref{tab:VGparams_combined}.  The pressure head distribution at different time steps is shown in Figure \ref{fig:testc3}. The pressure head distribution for the low order result overestimates the extent of wetting front due to the diffusion added in the low order method. The pure Galerkin method does not converge for this test problem. 

\subsection{Convergence Study Using Two-Dimensional Analytical Solution of RE}
\label{sec:tracy-analytical}
The analytical solutions derived by \citeA{https://doi.org/10.1029/2005WR004638} will be used to verify the convergence rate of the new FCT scheme. The Kirchhoff transformation along with Gardner relations is amenable to analytical solution techniques under specific Dirichlet boundary conditions using seperation of variables and Fourier transforms. The derived analytical solutions are for two dimensional and three dimensional cases. In this study, we will compare the transient two dimensional numerical solution with the transient analytical solution and determine the convergence rate.  

Let $\psi(\mathbf{x},t)$ denote pressure head, $\theta(\psi)$ water content, $K(\psi)$ the unsaturated hydraulic conductivity, and $\mathbf{e}_z$ the upward unit vector. The constitutive relationships that are used for the formulation of the analytical solutions are determined from the KSP relationships established by \citeA{gardnerSTEADYSTATESOLUTIONSUNSATURATED1958}  and \citeA{HydraulicConductivityUnsaturated1954}

\begin{equation}
    k_r= e^{\alpha \psi} \quad \quad
    \theta = \theta_r +(\theta_s - \theta_r) e^{\alpha \psi}.
\end{equation}

Define the Kirchhoff variable
\begin{linenomath*}
    \begin{equation}
        \bar h \;=\; e^{\alpha\psi} - e^{\alpha \psi_r},        
    \end{equation}
\end{linenomath*}

where $\alpha>0$ is a soil parameter and $\psi_r$ is initial pressure head. Written in terms of $\bar h$, the RE reduces to a linear convection–diffusion equation
\begin{equation}
\nabla^2 \bar h \;+\; \alpha\,\frac{\partial \bar h}{\partial z}
\;=\;
c\,\frac{\partial \bar h}{\partial t},
\qquad
c \;=\; \frac{\alpha(\theta_s-\theta_r)}{K_s}.
\label{eq:tracy-linear}
\end{equation}

We consider a rectangular domain $(x,z)\in[0,a]\times[0,L]$ with impermeable side
boundaries, a prescribed head at the top $z=L$, and a fixed reference head at the bottom
$z=0$. The initial condition is uniform dryness, $\psi(x,z,0)=\psi_r$, so that the transient response is driven entirely by the top boundary
forcing given in Table \ref{tab:tracy_parameters}.

\begin{table}[h]
\centering
\caption{Parameters and boundary conditions used in the steady-state and transient Tracy solution.}
\begin{tabular}{llc}
\hline
\textbf{Parameter} & \textbf{Description} & \textbf{Value} \\
\hline
$a$          & Domain width    & $10\ \mathrm{m}$ \\
$L$          & Domain height                          & $10\ \mathrm{m}$ \\
$\alpha$     & Gardner  parameter        & $0.164\ \mathrm{m}^{-1}$ \\
$\psi_r$     & Initial (reference) pressure head      & $-15.24\ \mathrm{m}$ \\
$K_s$        & Saturated hydraulic conductivity       & $2.04\ \mathrm{m\,d^{-1}}$ \\

\\
\hline
\textbf{Boundary Conditions} & & \\
Top ($z=L$)    & Dirichlet: 
$\displaystyle \psi_{\text{top}}(x)=\frac{1}{\alpha}
\ln\!\left[e^{\alpha\psi_r}
+ 0.5\,h_0\left(1-\cos\frac{2\pi x}{a}\right)\right]$ & -- \\
Bottom ($z=0$) & $-\psi_r$ & -- \\
Sides ($x=0,a$) & No Flow          & $\partial\psi/\partial z = 0$\\
\hline
\end{tabular}
\label{tab:tracy_parameters}
\end{table}

The steady-state solution is given by \eqref{eqn: steady_state_tracy}
\begin{linenomath*}
\begin{equation}
\label{eqn: steady_state_tracy}
\psi_{ss}(x,z)=
\frac{1}{\alpha}
\ln\!\left(
e^{\alpha\psi_r} + \bar{h}_{ss}(x,z)
\right),
\end{equation}    
\end{linenomath*}

\begin{linenomath*}
\begin{equation}
\bar{h}_{ss}(x,z)
=
\frac{\bar{h}_0}{2}\,
e^{\frac{\alpha}{2}(L-z)}
\left[
\frac{\sinh\!\left(\tfrac{\alpha}{2}z\right)}{\sinh\!\left(\tfrac{\alpha}{2}L\right)}
-
\cos\!\left(\frac{2\pi x}{a}\right)\,
\frac{\sinh(\beta_1 z)}{\sinh(\beta_1 L)}
\right]
\end{equation}    
\end{linenomath*}

where
\begin{linenomath*}
    \begin{equation}
        \bar{h}_0 = 1 - e^{\alpha \psi_r}, 
        \qquad
        \beta_1 = \sqrt{\left(\frac{\alpha}{2}\right)^2 + \left(\frac{2\pi}{a}\right)^2 }.        
    \end{equation}
\end{linenomath*}

The top boundary condition controls the infiltration intensity through the amplitude $\hat{h}_0$, while the depth dependence and lateral variation are encoded through the vertical exponential term and horizontal cosine term, respectively.  For the convergence study, parameters that are used for both analytical and the numerical setup are provided in Table \ref{tab:tracy_parameters}.

The transient solution is obtained by adding a decaying series to the steady-state solution. Let
\begin{equation}
\bar{\psi}(x,z,t) = \psi_{ss}(x,z) + \bar{\phi}(x,z,t).
\end{equation}
The vertical eigenvalues ($\lambda_k$) and the decay rates($\gamma_k$) for $k=1,2,\dots$ are provided below. 
\begin{equation}
\lambda_k = \frac{k\pi}{L},  \qquad \gamma_1 = \frac{1}{c}\left(\lambda_k^2 + \frac{\alpha^2}{4}\right), \qquad \gamma_2 = \frac{1}{c}\left(\left(\frac{2\pi}{a}\right)^2 + \frac{\alpha^2}{4}\right).
\end{equation}

\begin{figure}[ht]
    \centering
    \includegraphics[width=1.1\textwidth]{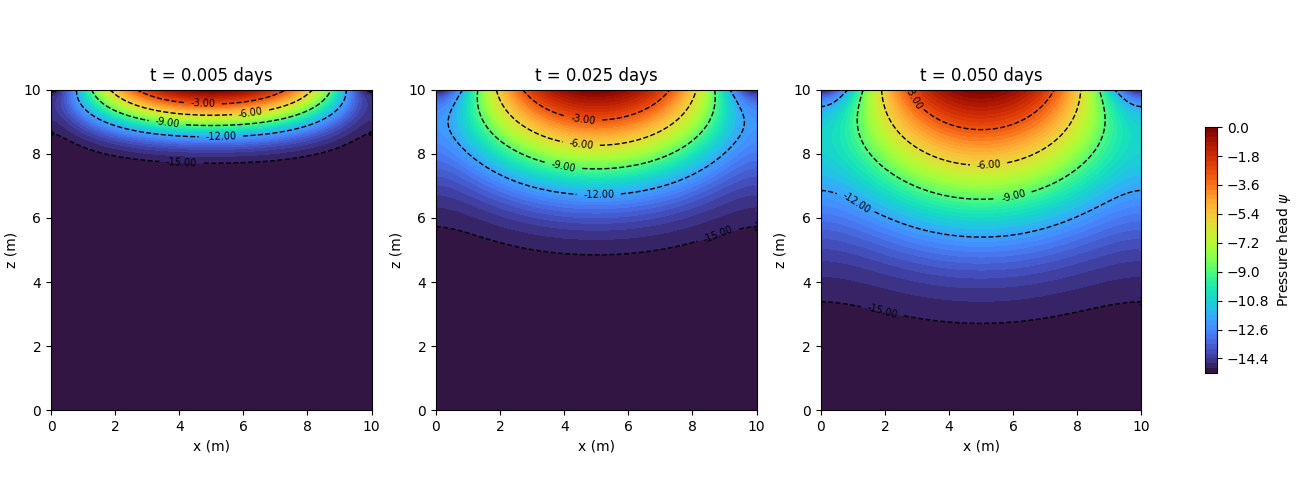}
    \caption{Transient solutions at different times based on the Tracy analytical solution.}
    \label{fig:Transient_analytical}
\end{figure}

The transient correction term and pressure head is
\begin{equation}
\psi(x,z,t)
=
\frac{1}{\alpha}
\ln\!\left[
e^{\alpha\psi_r}
+
\bar{h}_{ss}(x,z)
+
\bar{\phi}(x,z,t)
\right],
\label{eq:tracy-full-psi-unsteady}
\end{equation}
where, 

\begin{equation}
\bar{\phi}(x,z,t)
=
\frac{\bar{h}_0}{L\,c}\,
e^{\frac{\alpha}{2}(L-z)}
\sum_{k=1}^{\infty}
(-1)^k \lambda_k
\left[
\frac{1}{\gamma_1}
e^{-\gamma_1 t}\,\sin(\lambda_k z)
-
\frac{1}{\gamma_2}
e^{-\gamma_2 t}\,
\cos\left(\frac{2\pi x}{a}\right)
\sin(\lambda_k z)
\right].
\label{eq:tracy-transient}
\end{equation}
The parameters used for generating the transient analytical solution are provided in Table \ref{tab:tracy_parameters}.  The transient analytical solution at different time steps are shown in Figure \ref{fig:Transient_analytical}

The comparison between analytical and numerical solution at different refinement is shown in Table \ref{tab:transient_L2_convergence_comparison}. The table presents the transient spatial $L^2$ errors and convergence rates for the Galerkin and low Order schemes at several time instances. Across all mesh refinements, the Galerkin method consistently achieves approximately second-order accuracy, which shows the ability to capture the transient Tracy solution. On the other hand, The low order method produces significantly larger errors and exhibits first order convergence. The FCT solutions are substantially more accurate  and attain a second order convergence rate on the finer mesh. Importantly, the FCT scheme is fundamentally a predictor-corrector scheme for transient problems. While the objective of this work is a method for transient infiltration, particularly  for stormwater modeling, future work should consider methods that preserve bounds and perform well for steady-state variably-saturated flow. 

\begin{table}[h]
\centering
\begin{tabular}{c|cc|cc|cc}
\hline
\multirow{2}{*}{$h_e$} 
& \multicolumn{2}{c|}{\textbf{Galerkin }} 
& \multicolumn{2}{c|}{\textbf{Low Order}}
& \multicolumn{2}{c}{\textbf{FCT}} \\
\cline{2-7}
& $\|e\|_{L^2}$ & Order
& $\|e\|_{L^2}$ & Order
& $\|e\|_{L^2}$ & Order \\
\hline
\hline

\multicolumn{7}{c}{\textbf{$t = 3.75\times10^{-4}$ days }} \\
\hline
$0.1250$   & $4.743093$ & --
           & $1.491683$ & --
           & $1.918816$ & -- \\
$0.0625$   & $1.331997$ & $1.832$
           & $0.737184$ & $1.017$
           & $0.861955$ & $1.155$ \\
$0.03125$  & $0.301804$ & $2.142$
           & $0.378692$ & $0.961$
           & $0.193956$ & $2.152$ \\
$0.015625$ & $0.064863$ & $2.218$
           & $0.198489$ & $0.932$
           & $0.048290$ & $2.006$ \\
\hline\hline

\multicolumn{7}{c}{\textbf{$t = 4.375\times10^{-4}$ days }} \\
\hline
$0.1250$   & $4.789856$ & --
           & $1.440297$ & --
           & $1.848597$ & -- \\
$0.0625$   & $1.215065$ & $1.979$
           & $0.712814$ & $1.015$
           & $0.840755$ & $1.137$ \\
$0.03125$  & $0.273888$ & $2.149$
           & $0.366724$ & $0.959$
           & $0.198751$ & $2.081$ \\
$0.015625$ & $0.059041$ & $2.214$
           & $0.192198$ & $0.932$
           & $0.048044$ & $2.049$ \\
\hline\hline

\multicolumn{7}{c}{\textbf{$t = 5.0\times10^{-4}$ days }} \\
\hline
$0.1250$   & $4.784264$ & --
           & $1.397881$ & --
           & $1.790696$ & -- \\
$0.0625$   & $1.120624$ & $2.094$
           & $0.692613$ & $1.013$
           & $0.822607$ & $1.122$ \\
$0.03125$  & $0.251717$ & $2.154$
           & $0.356681$ & $0.957$
           & $0.203059$ & $2.018$ \\
$0.015625$ & $0.054404$ & $2.210$
           & $0.186918$ & $0.932$
           & $0.048256$ & $2.073$ \\
\hline
\end{tabular}
\caption{Transient convergence comparison (pressure head $\psi$) using the Tracy analytical solution.}
\label{tab:transient_L2_convergence_comparison}
\end{table}

\section{Stormwater Management Applications}
Stormwater infiltration systems, such as rain gardens and bioswales, are widely used for managing runoff by enhancing groundwater recharge and reducing surface water pollution. Due to the complexity of two-dimensional simulations, such systems are 
commonly designed using inaccurate models that rely on the assumption of 
1D saturated flow. Early infiltration systems were simulated as a one-dimensional horizontal flow based on the assumption that the flows would occur horizontally through the side walls of the trenches as the bottom will clog completely~\cite{mikkelsen1996infiltration}. On the other hand, many existing models like MUSIC(v3) stormwater model assume that the flow is in the 
vertical direction. Moreover, few models allow for quasi-2D conditions with vertical flow from the case and horizontal flow from the side walls~\cite{duchene1994modeling}. Most of these methods assume a constant infiltration rate or soil equations like Green Ampt model and Horton model~\cite{wong2006water, papa2005analysis}. Simulation of 2D infiltration models is required for accurate estimation of infiltration. However, very few 2D variably saturated models are used for the design of LID structures due to the complexity of the model and the unavailability of the data.  A horizontal one-dimensional model underestimated the flow from 24\% to 54\% compared to a 2D finite element model \cite{duchene1994modeling}. Furthermore, the assumption of one-dimensional vertical flow is even more conservative which underestimates the flow by 50\%. This reflects the importance of a two-dimensional model of infiltration in storm-water systems. In the study, two-dimensional simulations of rain gardens and bioswales are conducted for simulating the infiltration in the stormwater infiltration structures. 

Bioswales are designed to mimic the natural processes of infiltration. They are typically topped with native vegetation, which provides hydraulic resistance to stormwater flow. Although there are numerous studies on the treatment capacity of bioswales, there are not many investigations regarding the infiltration capabilities of bioswale. To date, most of the simulations are limited to macroscale water balances. Among them, SWMMM shows the best results for the design of bioswales that ensure optimal land use ~\cite{xu2017modelling}. In addition, DRAINMOD has been used effectively to assess the hydrologic response of bioswales to influent runoff \cite{brown2013calibration}.

\subsection{Rain Garden}
In this study, a two-dimensional numerical model of a rain garden was developed using the new FCT method. The numerical scheme incorporates adaptive time-stepping and is solved using Newton’s iteration with a tolerance criterion of $10^{-8}$. The time discretization was handled using Crank Nicolson method, and triangular elements were used to discretize the spatial domain. 

\begin{figure}[ht]
\centering
a)\includegraphics[width=0.30\textwidth]{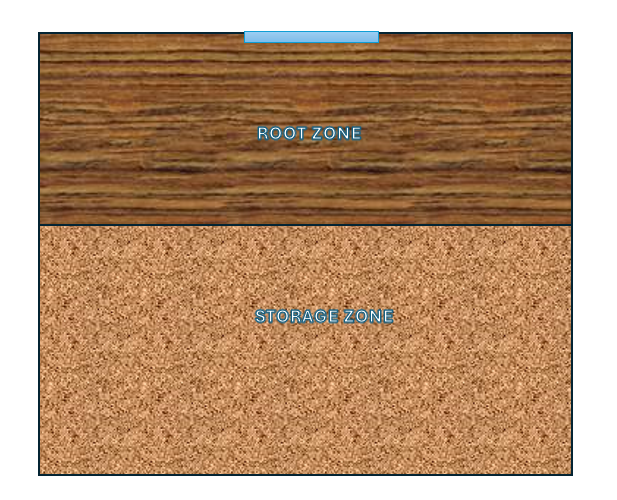}
b)\includegraphics [width=0.35\textwidth]{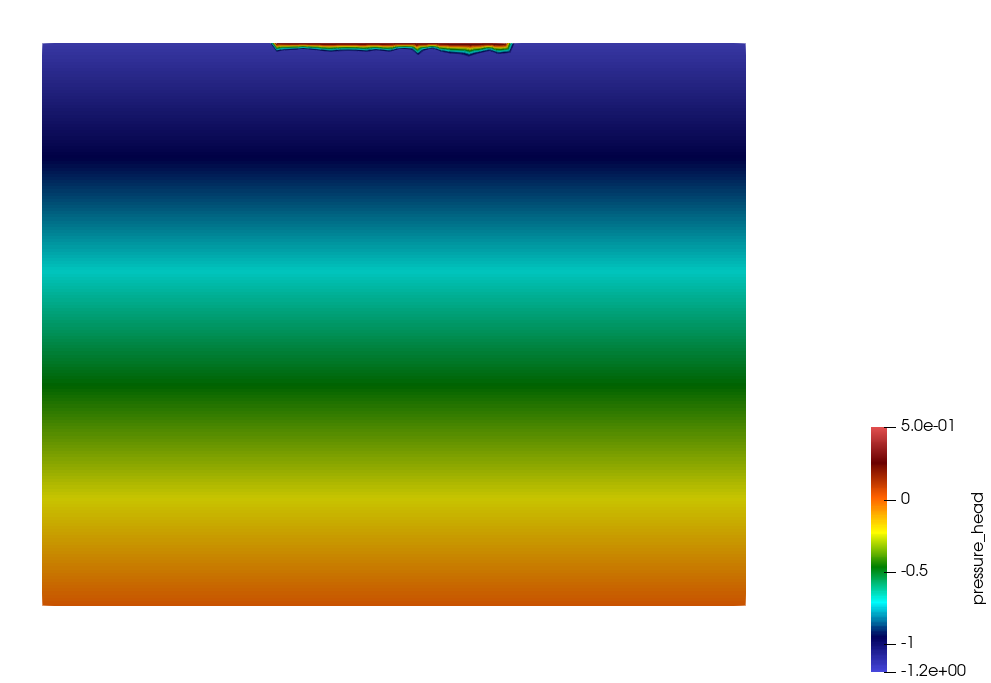}
\caption{Schematic diagram of the subsoil of rain garden, b) Initial  condition of rain garden.}
\label{fig: rg_schematic}
\end{figure}

Figure \ref{fig: rg_schematic} shows a two-dimensional sketch of the rain garden subsoil domain and the initial condition. The root zone of the rain garden consists of organic sand with a depth of 50cm and the storage zone consists of gravel or medium sand with a depth of 70 cm. The hydraulic properties for both layers were parameterized using the Mualem–van Genuchten model. The corresponding parameters are summarized in Table~\ref{tab:vg_params_rg}. The initial condition assumed a hydrostatic pressure distribution across the vertical profile. The top boundary was subject to a constant pressure head corresponding to a ponding depth of $\Psi_{\text{top}} = 0.5$ m. The sidewalls were treated with no-flow boundary conditions.
\begin{table}[ht]
    \centering
    \caption{Soil hydraulic properties for rain garden simulation.}
    \label{tab:vg_params_rg}
    \scalebox{0.95}{
    \begin{tabular}{lccccc}
    \hline
    \textbf{Layer} & $\alpha$ (m$^{-1}$) & $n$ & $\theta_r$ & $\theta_s$ & $K_s$ (m/day) \\
    \hline
    Root Zone        & 0.33   & 3.594   & 0.03  & 0.40  & 8.856 \\
    Storage Zone     & 0.32   & 2.146 & 0.10 & 0.37  & 19.944  \\
    \hline
    \end{tabular}
    }
\end{table}

\begin{figure}[ht] 
\centering
a)\includegraphics[width=0.45\textwidth]{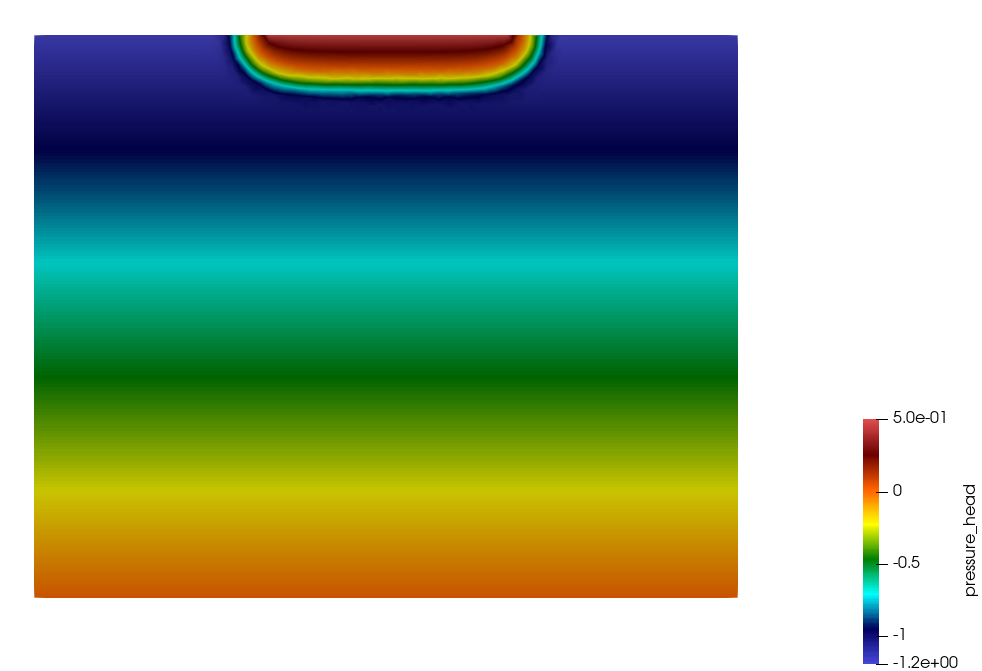}
b)\includegraphics[width=0.45\textwidth]{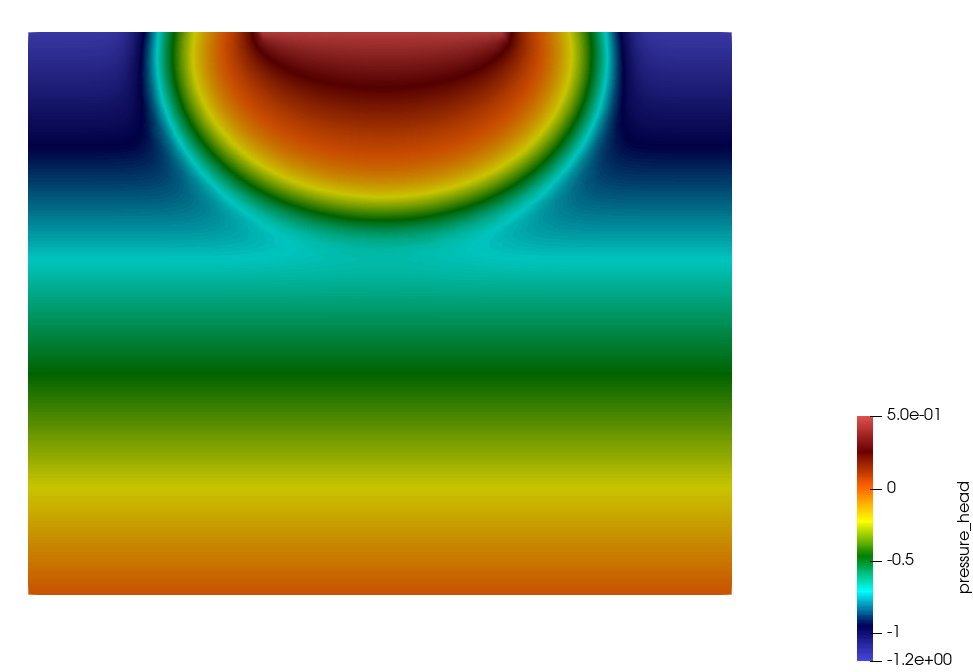}
c)\includegraphics[width=0.45\textwidth]{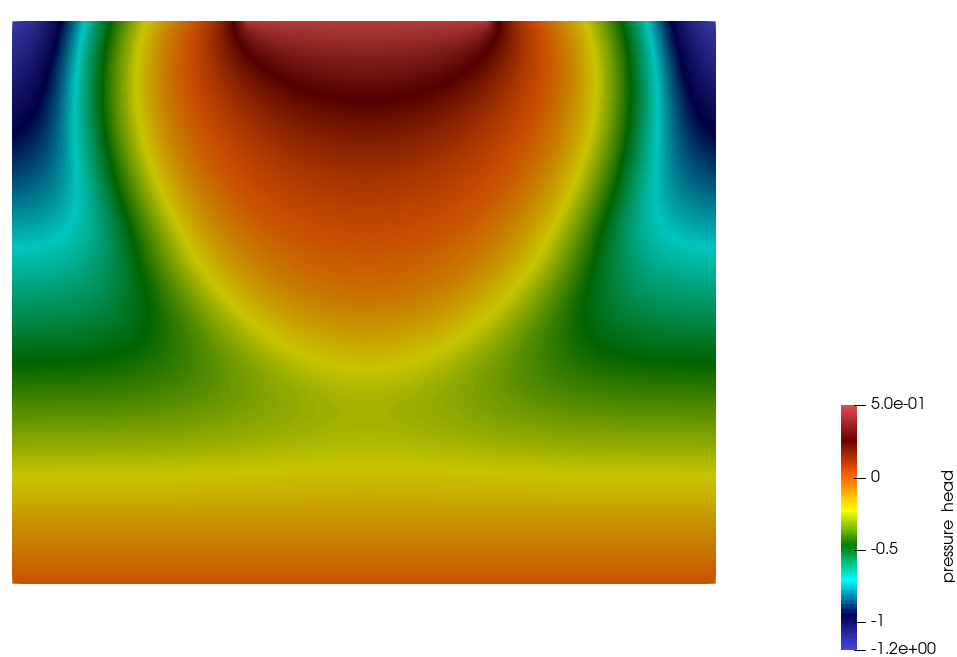}
d)\includegraphics[width=0.45\textwidth]{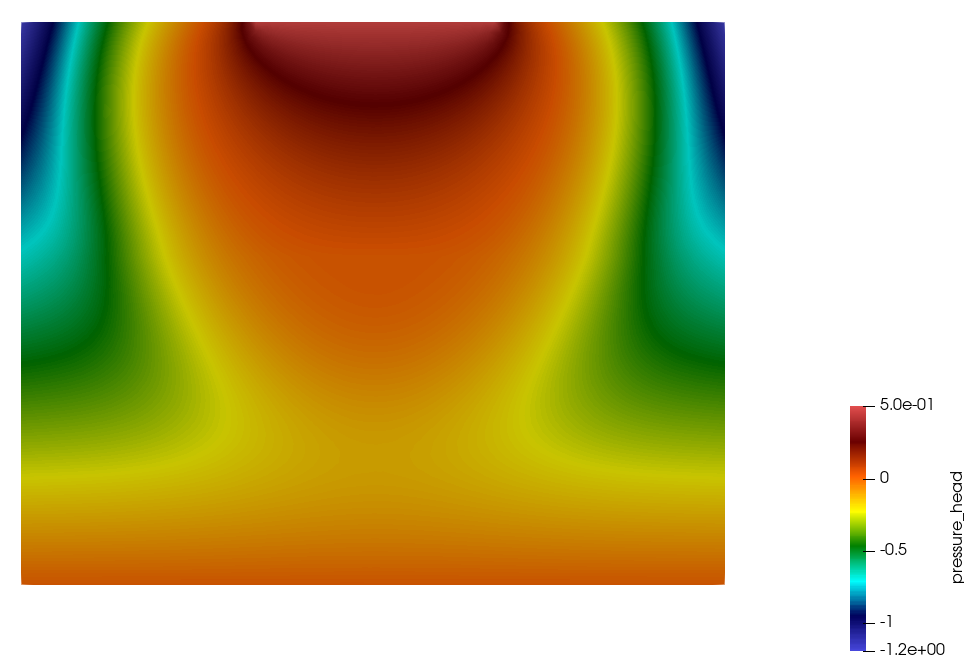}
\caption{ Simulation of infiltration in rain garden for a) 19.9 days, b) 199 
days, c) 497.512 days, d) 601.99 days.}
\label{fig: simulation_rg}
\end{figure}

Figure~\ref{fig: simulation_rg} illustrates the simulated pressure head distribution at an intermediate time during infiltration. The water infiltrates vertically from the surface ponding layer into the underlying soil as a result of gravitational and capillary forces. The subsoil becomes progressively saturated from top to bottom, whereas the sidewalls remain relatively dry due to slower lateral flow, demonstrating the dominance of vertical infiltration in this setting.

\subsection{ Bioswales} 
\begin{figure}[ht]
\centering
a)\includegraphics[width=0.4\textwidth]{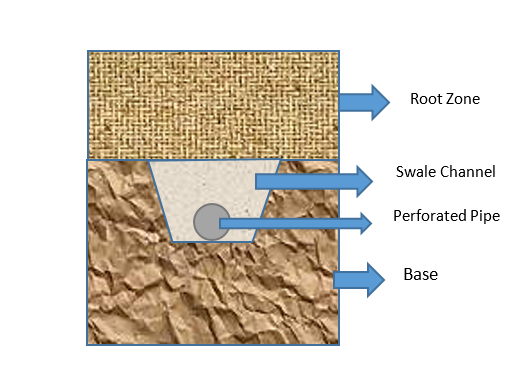}
b)\includegraphics[width=0.5\textwidth]{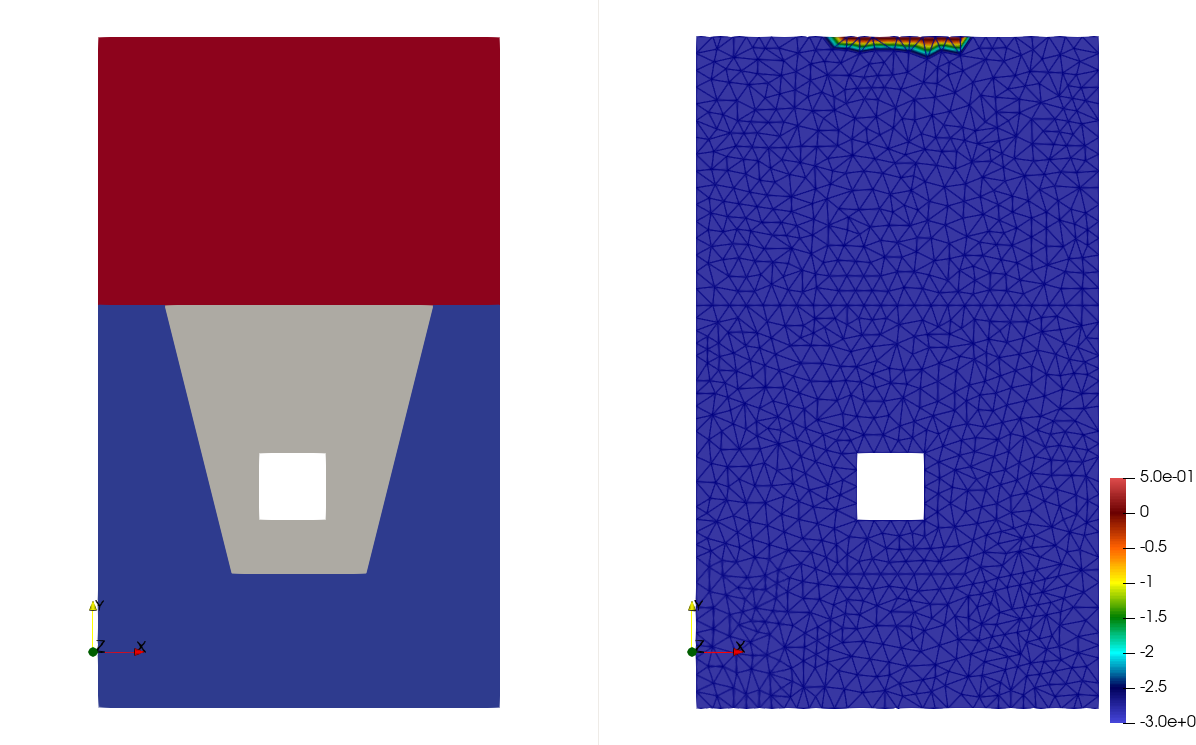}
\caption{a)Schematic Cross-section of the subsoil profile of the bioswale b) Material Type(left) and initial condition(right) of 2D bioswale.}
\label{fig: domain_initial_bio2d}
\end{figure}

In this study, the simulation domain consists of three distinct soil layers: the root zone, swale channel, and base zone, each characterized by varying hydraulic conductivities and retention properties (see Table~\ref{tab: bio_properties}). The presence of vegetation in the root zone increases the macropores, which results in improved hydraulic conductivity. On the other hand, the lower base layer is composed of finer materials with lower permeability. 
\begin{table}[!ht]
    \centering
    \caption{Soil properties of  bioswales.}
    \scalebox{0.8}{
    \begin{tabular}{llllll}
    \hline
        Zone & $K_s(m/d)$ & $\theta_r$ & $\theta_s$ & n & $\alpha(m^{-1})$ \\ 
        Root Zone & 5 & 0.05 & 0.4 & 2.4 & 8 \\ 
        Swale Channel & 7.128 & 0.045 & 0.43 & 2.68 & 14.5 \\ 
        Base & 1.06 & 0.065 & 0.41 & 1.89 & 7.5 \\ \hline
    \end{tabular}}
    \label{tab: bio_properties}
\end{table}
A schematic cross-section of the bioswale domain is shown in Figure~\ref{fig: domain_initial_bio2d}. The cross section of two-dimensional domain is 3x5m. The initial pressure head in the domain is uniformly set to $-3.00$~m, except in the ponded region at the surface which is located between $x = 1$~m and $x = 2$~m along the top boundary ($z = 5$~m). A constant pressure head of $\psi = 0.5$~m is applied at the ponded region and seepage boundary condition is imposed at the perforated pipe in the swale zone. In addition to that, no flow boundary is imposed to the other boundaries. The computational mesh for mesh diameter $h_e =$ 0.15m and the resulting initial pressure distribution are illustrated in Figure~\ref{fig: domain_initial_bio2d}.


\begin{figure}[ht]
\centering
a)\includegraphics[width=0.29\textwidth]{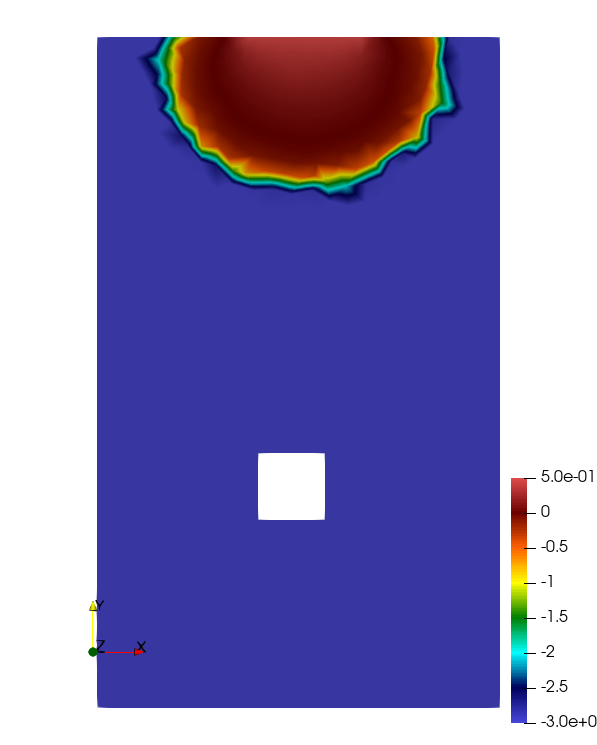}
b)\includegraphics[width=0.3\textwidth]{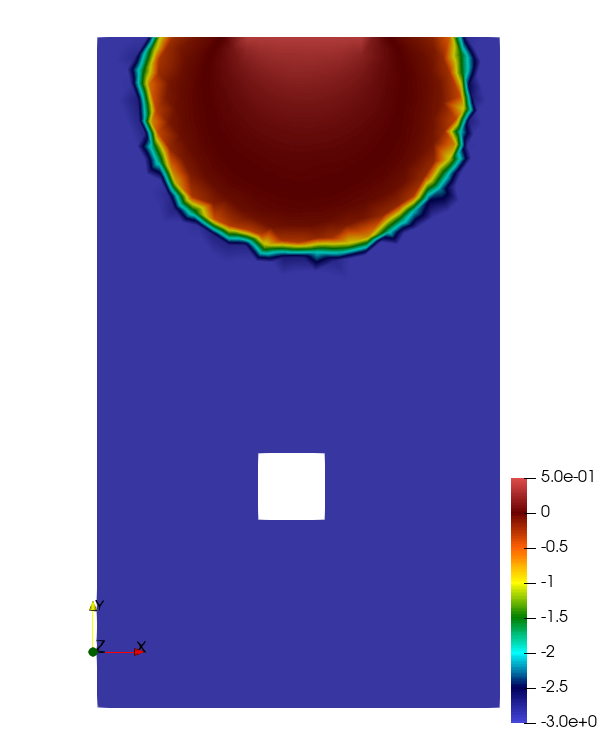}
c)\includegraphics[width=0.3\textwidth]{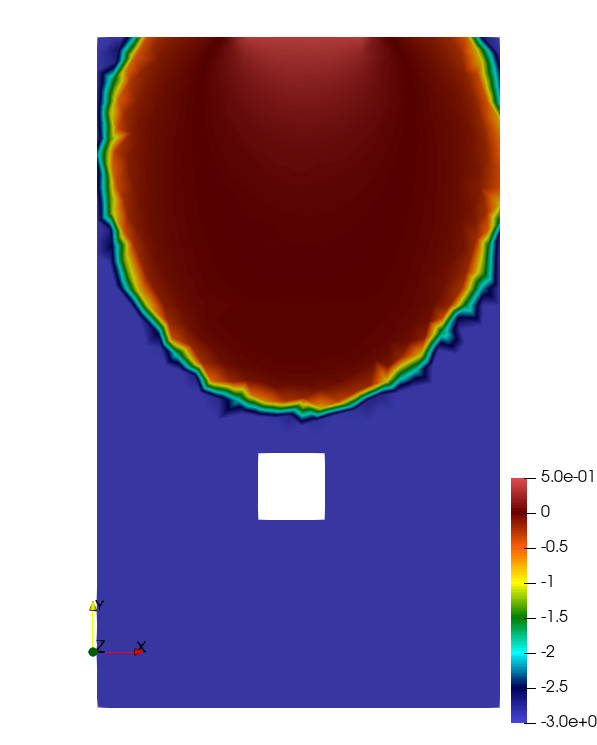}
d)\includegraphics[width=0.3\textwidth]{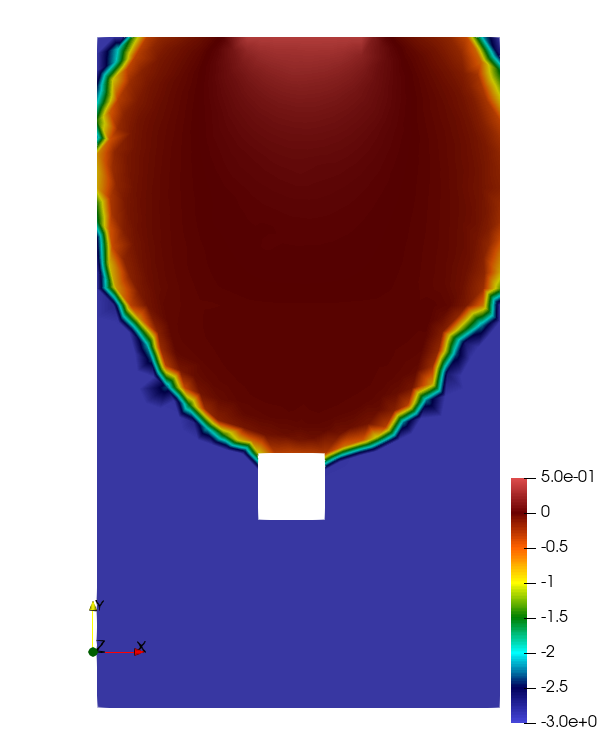}
e)\includegraphics[width=0.3\textwidth]{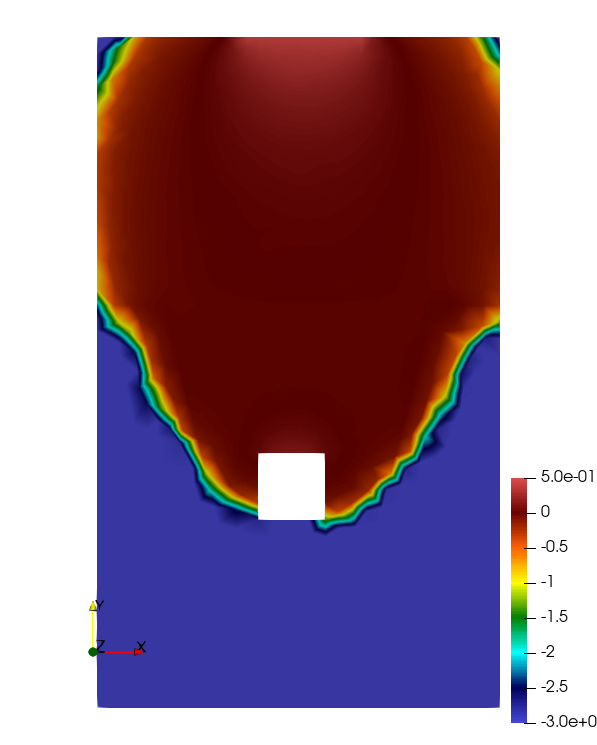}
f)\includegraphics[width=0.3\textwidth]{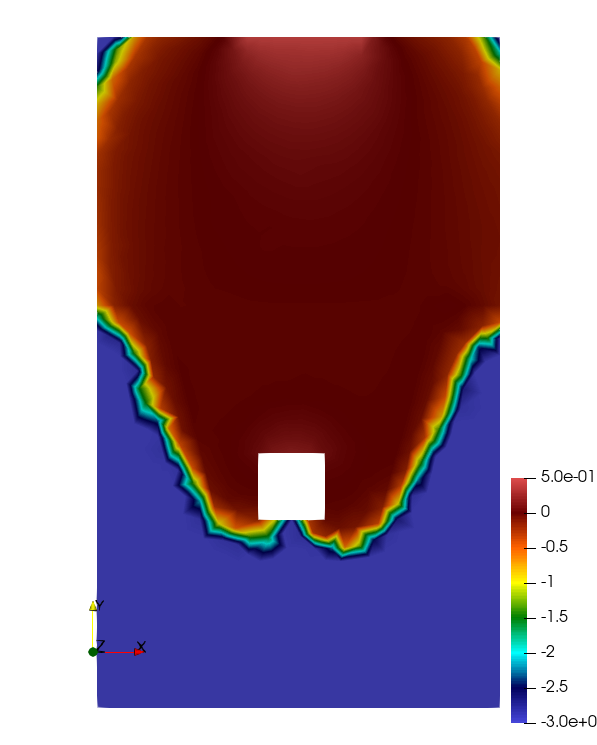}
\caption{Pressure head distribution in two-dimensional bioswale at a) 0.03 days b) 0.06 days c) 0.15 days d) 0.18 days e) 0.21 days f) 0.225 days.}
\label{fig: simbio2d}
\end{figure}

The simulations in Figure \ref{fig: simbio2d} provide a detailed temporal visualization of how the water front evolves and advances towards the perforated pipe embedded within the bioswale. Initially, the water infiltrates vertically down, following the gravitational pull and exhibits significant saturation across the domain. As the water front approaches the swale channel, it begins to curve below the perforated pipe. This behavior illustrates the interaction between the soil's permeability and the hydraulic gradients driving the flow. 

A mesh refinement study is also conducted to evaluate the accuracy of the proposed scheme. The high-resolution HYDRUS simulation with a mesh size of $h_e= 0.006m$ is used as a benchmark solution for the convergence study. Simulations are performed for $h_e = 0.3$, $0.15$, and $0.075$~m. The computed solutions are interpolated onto the HYDRUS mesh to evaluate relative $L_2$ errors. As shown in Table~\ref{tab:l2_error_convergence_multi}, the error decreases consistently with mesh refinement, demonstrating the effectiveness of the proposed method. At early times, convergence rates are greater than 2, which indicates high accuracy for smooth solution regimes. The benchmark solution is itself obtained from a high-resolution numerical simulation using HYDRUS; therefore, the convergence reflects consistency with a reference solution rather than comparison to an exact analytical solution.

\begin{table}[H]
\centering
\caption{Relative $L_2$ errors and convergence rates for different mesh sizes $h_e$.}
\label{tab:l2_error_convergence_multi}
\scalebox{0.95}{
\begin{tabular}{c|cc|cc}
\hline
\textbf{$h_e$} & \multicolumn{2}{c|}{$t = 0.06$} & \multicolumn{2}{c}{$t = 0.15$} \\
               & \textbf{Rel. Error} & \textbf{Rate} & \textbf{Rel. Error} & \textbf{Rate} \\
\hline
0.3   & 0.095806 & --    & 0.254473 & --    \\
0.15  & 0.014438 & 2.73  & 0.027596 & 3.20  \\
0.075 & 0.005989 & 1.27  & 0.008024 & 1.78  \\
\hline
\end{tabular}}
\end{table}

To further investigate the infiltration dynamics within stormwater management systems, a three-dimensional model of a bioswale is simulated. The domain extends 5 m in length. The simulation setup mirrors the 2D bioswale structure, including a layered soil, a vegetated root zone and a swale channel layer with coarse material as provided in Table~\ref{tab: bio_properties}. The initial condition of the pressure head is hydrostatic throughout the domain. The ponding is imposed at the top as a Dirichlet boundary condition. 

\begin{figure}[ht]
\centering
a)\includegraphics[width=0.45\textwidth]{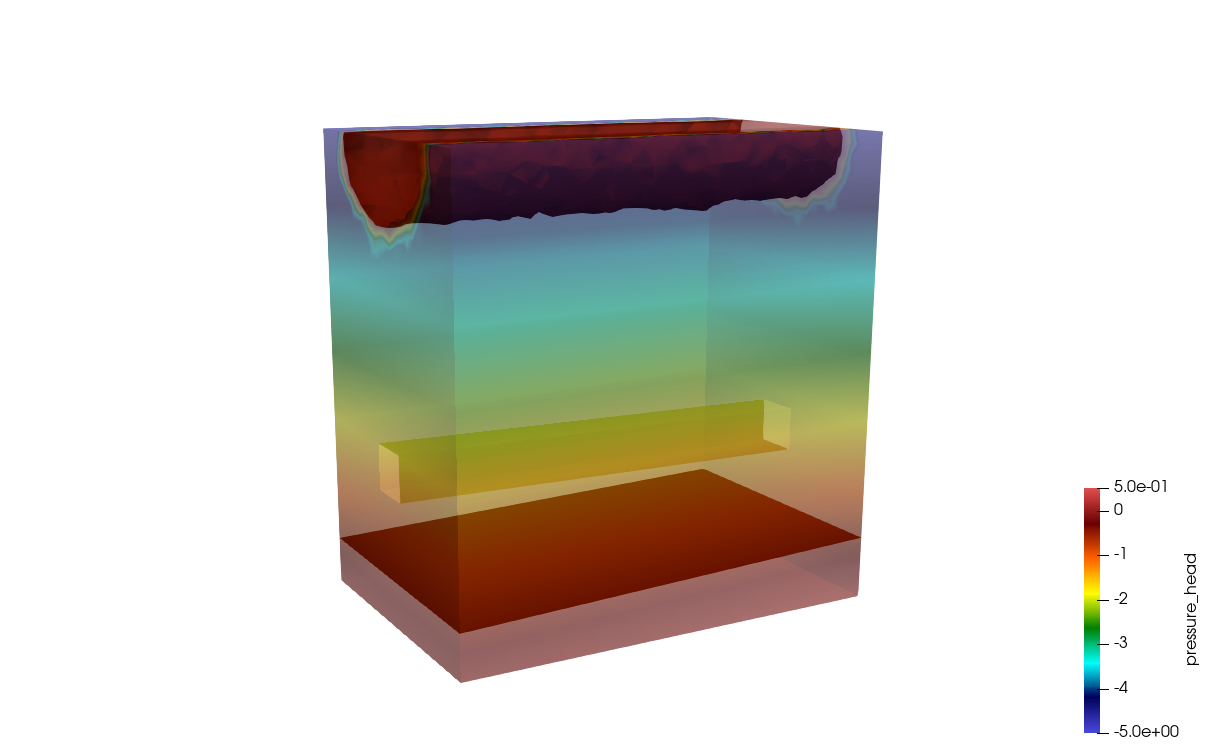}
b)\includegraphics[width=0.45\textwidth]{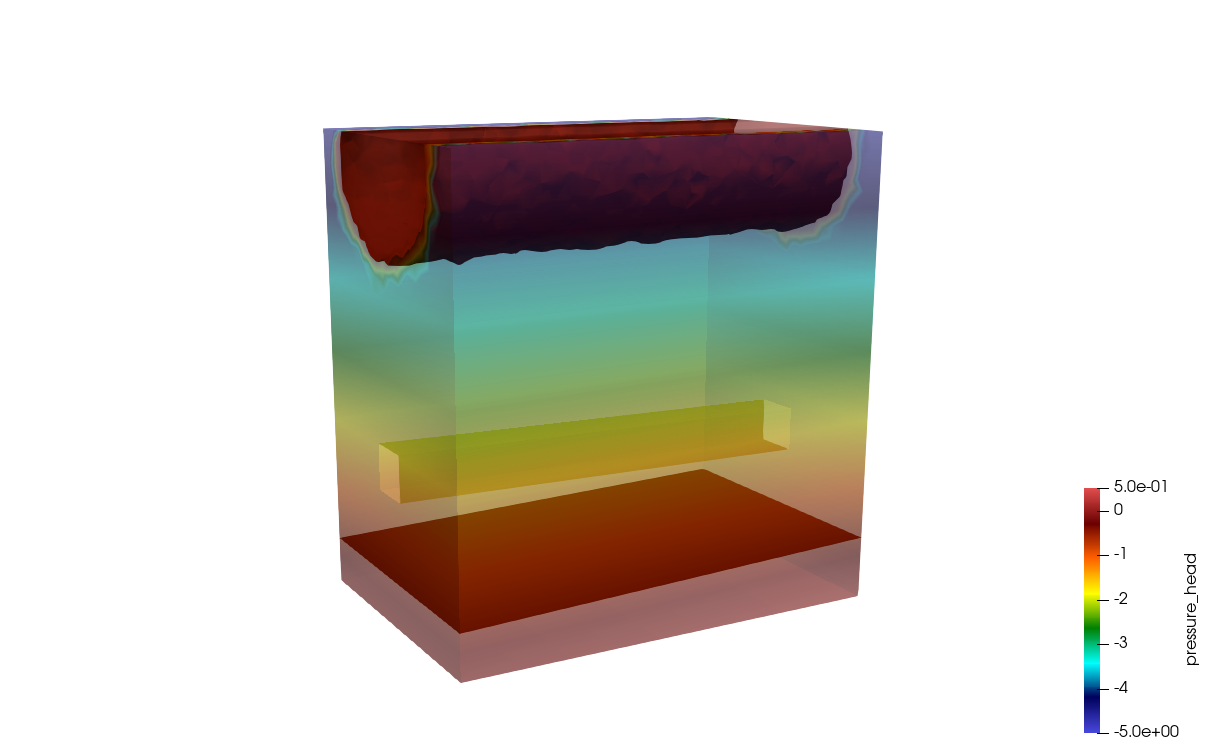}
c)\includegraphics[width=0.45\textwidth]{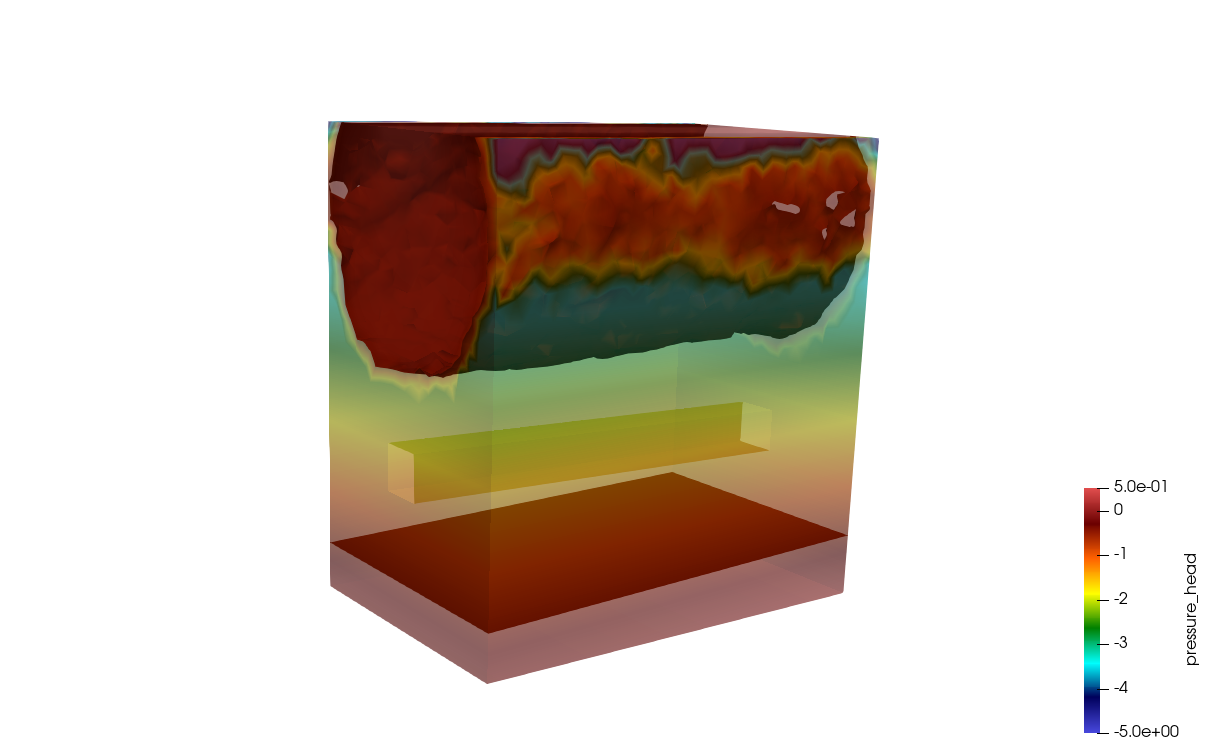}
d)\includegraphics[width=0.45\textwidth]{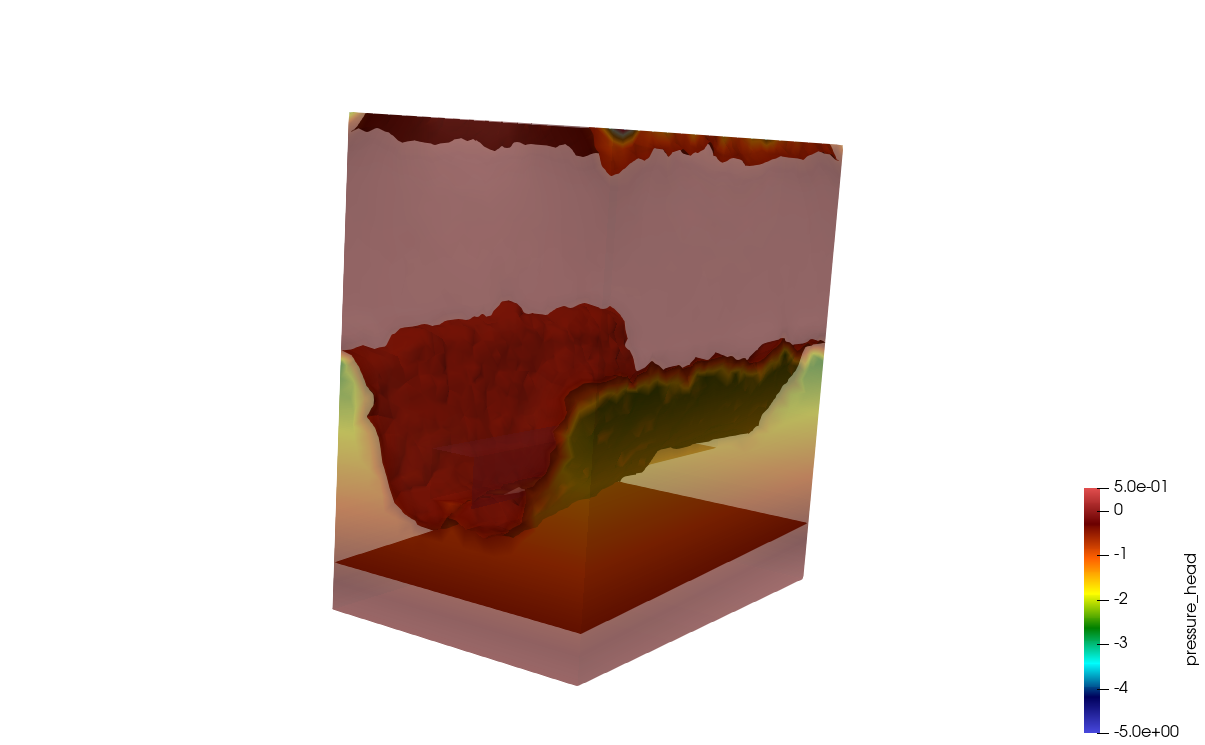}

\caption{Pressure head distribution in three-dimensional bioswale at a) 0.015 days b) 0.030 days c) 0.105 days d)0.21 days.}
\label{fig: simbio3d}
\end{figure}

Figure~\ref{fig: simbio3d} illustrates the temporal evolution of the pressure head distribution at several time steps. Similarities to the two-dimensional case are observed. The water changes its shape to advance to the swale channel under constant ponding at the top of the root zone. The results highlight the capability of the proposed model to simulate realistic three-dimensional flow pathways in engineered infiltration systems.

\section{Conclusion}

A novel version of FCT was presented for unsteady initial-boundary value problems for RE. The method has three components: 1) a fully implicit, bound-preserving, first-order predictor, based on well-known modifications of the piecewise-linear FEM, with existing rigorous analysis \cite{forsythMonotonicityConsiderationsSaturated1997,aithammououlhajNumericalAnalysisNonlinearly2018}. 2) An explicit FCT step to correct the moisture content while preserving local and global bounds.  And finally 3) a nodal Newton solve to recover the bound-preserving, second-order pressure. The properties of the new method were demonstrated quantitatively on several benchmarks and qualitatively on stormwater management applications. An important limitation is that the current scheme is not appropriate for steady-state calculations due to the predictor-corrector structure and the dependence of the limiting procedure on the time step size. On the other hand, the robustness of the scheme with respect to mesh quality and the ability to leverage existing FCT frameworks supporting third order and higher FEM, suggest promising future work to enable higher accuracy on coarser and/or adaptive meshes and design convergent methods for steady-state problems. 
\section{Open Research Section}
Source code for reproducing these results is  at \url{https://github.com/erdc/proteus}.

\acknowledgments

This work was funded by the United States Geological Survey through Louisiana Water Resources Research Institute grant GR-00010071.

%
%


\bibliography{RLn}

%
%
%
%
%

\end{document}